\documentclass[reqno, 11pt]{article}
\usepackage[all]{xy}
\usepackage{latexsym}
\usepackage{amsmath}
\usepackage[latin1]{inputenc}
\usepackage[T1]{fontenc}
\usepackage[francais]{babel}
\usepackage{mathdots}
\usepackage{amssymb,amscd,amsthm}
\usepackage[mathscr]{eucal}
\usepackage{bbm}
\usepackage{hyperref}
\newcommand{\F}{\mathbb F}

\newcommand{\N}{\mathbb N}

\newcommand{\Q}{\mathbb Q}
\newcommand{\R}{\mathbb R}
\newcommand{\Z}{\mathbb Z}

\newcommand{\rRep}{\mathrm{Rep}}
\newcommand{\rIrr}{\mathrm{Irr}}

\def\dim{\mathrm{dim}}
\def\Ker{\mathrm{Ker}}

\def\val{\mathrm{val}}

\def\Fil{\mathrm{Fil}}

\DeclareMathOperator{\Hom}{\mathrm Hom}

\DeclareMathOperator{\Frac}{\mathrm Frac}
\DeclareMathOperator{\GL}{\mathrm GL}

\DeclareMathOperator{\Gal}{\mathrm Gal}
\DeclareMathOperator{\id}{\mathrm Id}

\DeclareMathOperator{\Ind}{\mathrm {Ind}}
\DeclareMathOperator{\ind}{\mathrm {ind}}

\newcommand{\fixed}[1][V]{(B\otimes_{F}{#1})^{G}}

\newcommand{\xto}[1][]{\xrightarrow{#1}}

\newcommand{\simto}{%\stackrel{\sim}{\to}}%isomorphic to
\xto[\sim]} %isomorphic to;

\newcommand{\summ}{\sum\limits}

\newcommand{\bFp}{\overline{\F}_p}
\newcommand{\bQp}{\overline{\Q}_p}

\newcommand{\MOD}{\mathrm{MOD}}
\newcommand{\WD}{\mathrm{WD}}
\newcommand{\rec}{\mathrm{rec}}
\newcommand{\unit}{\mathrm{unit}}
\newcommand{\diag}{\mathrm{diag}}

\newcommand{\ide}{\mathbbm{1}}

\newtheorem{theorem}{Th\eb or\ed me}[section]
\newtheorem{lemma}[theorem]{Lemme}
\newtheorem{cor}[theorem]{Corollaire}
\newtheorem{prop}[theorem]{Proposition}

\theoremstyle{definition}
\newtheorem{defn}[theorem]{D\eb finition}
\newtheorem{conj}[theorem]{Conjecture}

\newtheorem{rem}[theorem]{Remarque}
\newtheorem{examples}[theorem]{Exemples}

\newcommand{\ra}{\rightarrow}

\newcommand{\ad}{\`{a}}

\newcommand{\eb}{\'{e}}
\newcommand{\ec}{\^{e}}
\newcommand{\ed}{\`{e}}

\newcommand{\ud}{\`{u}}

\begin{document}

\title{Normes invariantes et existence de filtrations admissibles}
\author{Yongquan Hu} %\\
\date{}
 \maketitle

\vspace{-10mm}

\hspace{5cm}\hrulefill\hspace{5.5cm} \vspace{5mm}

\textbf{Résumé} -- Dans \cite{BS},  est formulée une conjecture
sur l'équivalence entre l'existence de normes invariantes sur
certaines représentations localement algébriques de $\GL_{d+1}(L)$
et l'existence de certaines représentations de de Rham de
$\Gal(\bQp/L)$, où $L$ est une extension finie de $\Q_p$. Dans cet
article, on montre le sens ``facile'' de cette conjecture:
l'existence de normes invariantes entraîne l'existence de
filtrations admissibles.
\\

\textbf{Abstract} -- In \cite{BS}, is formulated a conjecture on the
equivalence of the existence of invariant norms on certain locally
algebraic representations of $\GL_{d+1}(L)$ and the existence of
certain de Rham representations of $\Gal(\bQp/L)$, where $L$ is a
finite extension of $\Q_p$. In this paper, we prove the ``easy''
direction of the conjecture\hspace{-1.5mm}: the existence of
invariant norms implies the existence of admissible filtrations.

%{\def\contentsname{Table des matières}%
\tableofcontents
\thanks{\noindent \hrulefill \\  \quad\textit{Mathematics Subject Classification (1991)}: 05A20,
17B10,
 22E50}

\section{Introduction}
L'objet de cet article est une conjecture proposée dans \cite{BS}.
Pour l'énoncer, on fixe d'abord quelques notations. Soit $p$ un
nombre premier. On fixe deux extensions finies  $L$ et $K$ de $\Q_p$
contenues dans $\overline{\Q}_p$ telles que
$[L:\Q_p]=|\Hom_{\Q_p}(L,K)|$, et aussi une  extension finie
galoisienne $L'$ de $L$ telle que
$[L_0':\Q_p]=|\Hom_{\Q_p}(L_0',K)|$ où $L_0'$ désigne le
sous-corps non-ramifié maximal de $L'$. On note
$W(\overline{\Q}_p/L)$ (resp. $W(\bQp/L')$) le groupe
de Weil de $L$ (resp. $L'$). %Notons $k_{L}$ le corps ré siduel de
%$L'$ et $|k_{L'}|=p^{f'}$.

 On considère les
deux catégories suivantes (cf. \S2.2):

(i) $\WD_{L'/L}$: la catégorie des $K$-représentations $(r,N,V)$
de dimension finie du groupe de Weil-Deligne de $L$ telles que la
restriction de $r$ à $W(\bQp/L')$ est non-ramifiée;

(ii) $\MOD_{L'/L}$: la catégorie des $(\varphi,N)$-modules étales
sur $L_0'\otimes_{\Q_p}K$ munis d'une action de $\Gal(L'/L)$
vérifiant quelques conditions usuelles.

%Ces deux catégories sont en fait équivalentes. Plus pré
%cisé ment,
On sait qu'il existe deux foncteurs (non canoniques)
%$\WD:\MOD_{L'/L}\ra \WD_{L'/L}$
% et $\MOD:\WD_{L'/L}\ra \MOD_{L'/L}$ tels qu'ils induite une équivalence entre ces
%deux catégories (cf. \S2.2).
\[\WD:\MOD_{L'/L}\ra \WD_{L'/L},\ \ \ \MOD:\WD_{L'/L}\ra \MOD_{L'/L}\]
qui sont inverses l'un de l'autre (cf. \S2.2).

Rappelons que (cf. \cite{Fo1}) pour un objet
$(\varphi,N,\Gal(L'/L),D)$ de $\MOD_{L'/L}$ on peut définir un
entier $t_N(D)$ et, si de plus on se donne d'une filtration
décroissante exhaustive séparée sur $D_{L'}=L'\otimes_{L_0'}D$ qui
est stable sous l'action de $\Gal(L'/L)$, on a un autre entier
$t_H(D_{L'})$. Par définition cette filtration est dite
\emph{admissible} si $t_H(D_{L'})=t_N(D)$ et si $t_H(D'_{L'})\leq
t_{N}(D')$ pour tout sous-objet $D'$ de $D$ où $D'$ est muni de la
filtration induite. Pour tout $\sigma:L\hookrightarrow K$, on pose
\[D_{L',\sigma}=D_{L'}\otimes_{L'\otimes_{\Q_p}K}(L'\otimes_{L,\sigma}K).\]
On a un isomorphisme naturel
$D_{L'}\simeq\prod_{\sigma:L\hookrightarrow K}D_{L',\sigma}$.

Si $(r,N,V)\in \WD_{L'/L}$, on note $(r,N,V)^{\textrm{ss}}\in
\WD_{L'/L}$ sa $F$-semisimplification (cf. \cite{De}, $\S$8.5).

On fixe un entier $d\geq 1$ et:

(i) un objet $(r,N,V)$ de $\WD_{L'/L}$ de dimension $d+1$ et tel que
$r$ est semi-simple;

(ii) pour tout $\sigma:L\hookrightarrow K$, un ensemble de $d+1$
entiers $i_{1,\sigma}<\cdots<i_{d+1,\sigma}$.

Posons $a_{j,\sigma}=-i_{d+2-j,\sigma}-(j-1)$ pour tout $\sigma$
et $j\in \{1,\dots,d+1\}$. D'une part, \`{a} $(r,N,V)$ on associe
une représentation lisse $\pi$ de $\GL_{d+1}(L)$ par la
correspondance de Langlands locale modifiée et un
$(\varphi,N)$-module sur $L_0'\otimes_{\Q_p}K$ muni d'une action
de $\Gal(L'/L)$ (cf. $\S$2.2). D'autre part, à
$\{a_{j,\sigma}\}_{j,\sigma}$ on associe une représentation
$\Q_p$-rationnelle sur $K$ de $\GL_{d+1}(L)$, notée $\rho$.

\begin{conj}(\cite{BS}, Conjecture 4.3)
Les deux conditions suivantes sont équivalentes:

 (i)
$\rho\otimes\pi$ admet une norme invariante;

(ii) Il existe un objet $(\varphi,N,\Gal(L'/L),D)$  dans
$\MOD_{L'/L}$ tel que:
\[WD(\varphi,N,\Gal(L'/L),D)^{\mathrm{ss}}=(r,N,V),\]
et une filtration admissible $(\Fil^i D_{L',\sigma})_{i,\sigma}$
stable par $\Gal(L'/L)$ sur $D_{L'}$ telle que
\[\Fil^i D_{L',\sigma}/\Fil^{i+1}D_{L',\sigma}\neq 0\Leftrightarrow i\in\{i_{1,\sigma},\dots,i_{d+1,\sigma}\}.\]
\end{conj}
%\begin{rem}
 %Soit $K'$ une extension finie galoisienne de $K$ et on pose
%$\rho\otimes_K{K'}$, $\pi\otimes_KK'$, $D\otimes_{K}K'$ les
%extensions des scalaires de $K$ à\ $K'$. Alors
%
%(1) La condition (i) est vraie pour $K$ si et seulement si (i)  est
%vraie pour $K'$.
%
%(2) La condition (ii) est varie pour $K$ si et seulement si (ii) est
%vraie pour $K'$ et la filtration sur $D_{L'}\otimes_K{K'}$ est
%stable par $\Gal(K'/K)$.
%\end{rem}

On renvoie le lecteur à \cite{BS}, $\S$1 pour l'introduction de la
conjecture. Certains résultats sur cette conjecture sont connus:
\begin{itemize}
\item[--] si $(r,N,V)$ est absolument irréductible (alors $N=0$), la conjecture est
vraie (\cite{BS}, Theorem 5.2);

\item[--] si $(r,N,V)$ est absolument indécomposable, alors la condition (ii) est vraie si et seulement si  le
caractère central de $\rho\otimes\pi$ est unitaire; en particulier,
on a (i) entraîne (ii) (\cite{BS}, Proposition 5.3);

\item[--] si $(r,N,V)$ est tel que $N=0$ et $r$ est non-ramifiée
et scindée sur $K$, alors (i) entraîne (ii) (\cite{BS}, Theorem
5.6);

\item[--] si $L=L'=\Q_p$ et $d=1$, la conjecture est vraie (\cite{BB}, \cite{Co}).
\end{itemize}

Dans cet article, nous montrons le théorème suivant:
\begin{theorem}
Supposons $K$ suffisamment gros tel que l'on a
\[(r,N,V)=\bigoplus_{i=1}^s(r_i,N_i,V_i)\]
avec les $(r_i,N_i,V_i)$ absolument indécomposables de dimension
$d_i$. On considère les quatre conditions suivantes:

(i) $\rho\otimes\pi$ admet une norme invariante.

(ii) Il existe un objet $(\varphi,N,\Gal(L'/L),D)$ dans
$\MOD_{L'/L}$ tel que:
\[WD(\varphi,N,\Gal(L'/L),D)^{\mathrm{ss}}=(r,N,V),\]
et une filtration admissible $(\Fil^i D_{L',\sigma})_{i,\sigma}$
 stable par $\Gal(L'/L)\times\Gal(K/\Q_p)$ sur $D_{L'}$ telle que
\[\Fil^i D_{L',\sigma}/\Fil^{i+1}D_{L',\sigma}\neq 0\Leftrightarrow i\in\{i_{1,\sigma},\dots,i_{d+1,\sigma}\}.\]

(iii) Posons $(\varphi_i,N_i,\Gal(L'/L),D_i)=\MOD(r_i,N_i,V_i)$.
Pour une permutation convenable $\nu$ de $\{1,\dots, s\}$ (voir
\S2.4), les inégalités suivantes sont vérifiées:
\[[K:L]\summ_{j=1}^{d_{\nu(1)}}\summ_{\sigma}i_{j,\sigma}\leq t_N(D_{\nu(1)}),\]
\[\vdots\]
\[[K:L]\summ_{j=1}^{d_{\nu(1)}+\cdots+d_{\nu(s-1)}}\summ_{\sigma}i_{j,\sigma}\leq\summ_{i=1}^{s-1}t_N(D_{\nu(i)}),\]
\[[K:L]\summ_{j=1}^{d+1}\summ_{\sigma}i_{j,\sigma}=\summ_{i=1}^{s}t_N(D_{\nu(i)})=t_N(D).\]

(iv)   Le caractère central de $\rho\otimes\pi$ est unitaire et la
condition d'Emerton est satisfaite (voir
$\S$\ref{section-Emerton}).

 Alors, on a les implications et équivalences:
 \[(i)\Rightarrow(ii)\Leftrightarrow(iii)\Leftrightarrow(iv).\]
\end{theorem}
\begin{rem}
L'existence de la permutation $\nu$ dans (iii) et l'équivalence
entre (ii) et (iii) sont la réponse de \cite{BS}, Remark 5.7.
\end{rem}

Comme conséquence, on obtient (cf. corollaire \ref{cor-main})
\begin{cor}
Avec les notations comme dans la conjecture 1.1, on a l'implication%(i) entraîe (ii).
\[(i)\Rightarrow(ii).\]
\end{cor}

\hspace{6mm}

 Voici le plan de l'article:

 On commence  par rappeler certaines des constructions et
notations dans \cite{BS} pour pouvoir énoncer la conjecture
précédente et le théorème principal ($\S$2). Les autres
chapitres sont consacrés à prouver notre théorème. Plus
précisément,
%on trouve d'abord une extension finie galoisienne $K'$
%de $K$ qui est suffisamment grande telle qu'on a
%\[(r,N,V)=\bigoplus_{i=1}^s(r_i,N_i,V_i)\]
%avec les $r_i$ absolument indécomposable, et
on prouve
\[(i)\Rightarrow(iv)\Leftrightarrow(iii)\Leftrightarrow(ii).\]
L'implication $(i)\Rightarrow(iv)$ est en fait un lemme d'Emerton
(cf. lemme \ref{lemme-Emerton}) et $(ii)\Rightarrow(iii)$ se déduit
de la définition d'admissibilité. Au $\S$3, on montre
$(iii)\Leftrightarrow (iv)$ en exprimant la condition d'Emerton
explicitement par la théorie de Bernstein-Zelevinsky.

Au \S4, on montre $(iii)\Rightarrow (ii)$.
% sous l'hypothè se que $K$ est suffisamment gros.
Cette partie, bien que constituée d'algèbre (semi-)%\linebreak
linéaire, est la plus technique de l'article. L'idée de la
démonstration est comme suit:  Après l'étude de quelques
exemples, on définit d'abord l'objet $(\varphi,N,\Gal(L'/L),D)$
(essentiellement l'opérateur $\varphi$).
%On étudie la structure de $D$ et introduit
%la notion \emph{bonté}\ d'un sous-objet $D'$
% de $D$ dans $\S$4.2. Remarquons que l'ensemble des bons sous-objets est de cardinal fini.
Au \S4.2 on introduit un certain ensemble fini de sous-objets de
$D$, noté $S$, et
% Dans le $\S$4.3
on munit $D_{L'}$ d'une filtration convenable (\S4.3) telle que pour
tout
 $D'\in S$ on a
\[t_H(D'_{L'})= \summ_{j=1}^{\textrm{rg}
D'}\summ_{\sigma}i_{j,\sigma},\] et par conséquent
$t_H(D'_{L'})\leq t_N(D')$. On montre que cette filtration est
\emph{admissible}. Pour cela, soit $D'$ un  sous-objet quelconque de
$D$, on lui associe une suite de sous-objets  $E_i\in S$
\[0=E_0\subsetneq E_1\subsetneq\cdots\subsetneq E_m\subsetneq E_{m+1}=D,\]
qui est la plus \emph{proche} de $D'$ au sens où  si $L\in S$ est
de même rang que $E_i$, alors
\[\textrm{rg}_{L_0'\otimes_{\Q_p}K} (D'\cap L)\leq \textrm{rg}_{L_0'\otimes_{\Q_p}K}(D'\cap E_i).\]
%De plus, cette propriété détermine le rang de $E_i$ sur
%$L_0'\otimes_{\Q_p}K$.
 D'autre part, si on pose $c_i=\textrm{rg}
(D'\cap E_{i})-\textrm{rg} (D'\cap E_{i-1})$ (sur
$L_0'\otimes_{\Q_p}K$) et
\[\Omega=\{j\in\N|\ \exists i \textrm{ tel que }\textrm{rg} E_i-c_i+1\leq j\leq \textrm{rg} E_i \},\]
 alors on a $t_H(D'_{L'})\leq
[K:L]\summ_{j\in\Omega}\summ_{\sigma}i_{j,\sigma}$. Enfin on compare
$[K:L]\summ_{j\in\Omega}\summ_{\sigma}i_{j,\sigma}$ et $t_N(D')$ en
utilisant un lemme combinatoire, et on en déduit que
$t_{H}(D'_{L'})\leq t_N(D')$ (cf. le corollaire
\ref{cor-admissible}).
\\

\textbf{Remerciement.} Ce travail s'est accompli sous la direction
de C. Breuil. Je le remercie pour avoir partagé avec moi ses idées
et ses connaissances et pour toutes ses remarques.

\section{Rappels et notations}
\subsection{Notations générales}
Soit $p$ un nombre premier. On fixe une clôture algébrique
$\overline{\Q}_p$ de $\Q_p$, et aussi deux extensions finies  $L$ et
$K$ de $\Q_p$ contenues dans $\overline{\Q}_p$ telles que
$[L:\Q_p]=|\Hom_{\Q_p}(L,K)|$. On note $\Gal(\overline{\Q}_p/L)$ le
groupe de Galois de $L$ et $W(\overline{\Q}_p/L)$ son groupe de Weil
(qui est dense dans $\Gal(\overline{\Q}_p/L)$), et
$\rec:W(\overline{\Q}_p/L)^{\mathrm{ab}}\simto L^{\times}$
l'isomorphisme  de réciprocité de telle sorte que les Frobenius
arithmétiques s'envoient sur les inverses des uniformisantes.
Posons $k_L$ le corps résiduel de $L$ et $q=p^f=|k_L|$. On note
$L_0=\Frac(W(\F_q))$ le sous-corps non-ramifié maximal de $L$ et
$\varphi_0$ le Frobenius sur $L_0$. On note $\val_{p}$ la valuation
$p$-adique sur $\overline{\Q}_p$ normalisée par $\val_{p}(p)=1$ et
on pose $|x|_p=p^{-\val_p(x)}$; de même on pose
$\val_L(x)=e\val_{p}(x)$ et  $|x|_L=q^{-\val_L(x)}$, où\
$e=[L:\Q_p]/f$ est l'indice de ramification de $L$ sur $\Q_p$.
%Si $w\in
%W(\overline{\Q}_p/L)$ et $\bar{w}$ son image dans
%$\Gal(\overline{\F}_p/\F_p)$, on note $\alpha(w)$ l'entier défini
%par $\overline{w}(x)=x^{p^{\alpha(w)}}$ pour $x\in\overline{\F}_p$.
%Remarquons que $\alpha(w)\in f\Z$.

Fixons $L'$ une extension finie galoisienne de $L$. On définit
$L_0'$, $k_{L'}$, $f'$, $\varphi_0'$, $\val_{L'}$ comme ci-dessus.
On suppose que $[L_0':\Q_p]=|\Hom_{\Q_p}(L_0',K)|$.

%On pose $G=\GL_{d+1}(L)$, et si $n\in\N$, on pose $G_n=GL_n(L)$. On
%note $\rRep G_n$ la catégorie des $\bQp$-représentations
%lisses de $G_n$ et $\rIrr G_n$ la sous-catégorie formée des
%représentations irréductibles.
%Posons $G=\GL_n(L)$, $\mathbf{G}=\mathbf{GL_n}$ et
%\[\tilde{\mathbf{G}}=(\mathrm{Res_{L/\Q_p}}\mathbf{G})_K\]
%le groupe réductif connexe changement de base de la restriction
%à\ la Weil. Si $\sigma:L\hookrightarrow K$ est un plongement, on
%note $\mathbf{G}_{\sigma}$ le changement de base de $\mathbf{G}$ par
%$\sigma$. On sait que
%\[\tilde{\mathbf{G}}=\prod_{\sigma:L\hookrightarrow K}\mathbf{G}_{\sigma}.\]

Si $\mathbb{G}$ est un groupe réductif défini sur $\Q_p$ et si
$G=\mathbb{G}(L)$ est le groupe des $L$-points rationnels de
$\mathbb{G}$, on note $\rRep G$ la catégorie des
$\bQp$-représentations lisses de $G$ et $\rIrr G$ la
sous-catégorie formée des représentations irréductibles. Soit
$P$ un sous-groupe parabolique de $G$ avec $M$ son quotient de Levi
et $N$ son radical unipotent. On définit les foncteurs suivants:
\[\ind_{P}^{G}: \ \ \Ind_{P}^{G}:\rRep M\ra\rRep G,\]
\[   r_{P}^{G}:\rRep G\ra \rRep M.\]

(a) Soit $(\sigma,W)\in\rRep G$, notons $\ind_{P}^G\sigma$ l'espace
des fonctions $f:G\ra W$ telles que
\begin{itemize}
\item[--] $f(nmg)=\sigma(m)f(g)$, si $n\in N$, $m\in M$ et $g\in G$,

\item[--] $f$ est invariante à droite par un sous-groupe ouvert
de $G$.
\end{itemize}
Le groupe $G$ opère par translation à droite et on obtient une
représentation lisse de $M$. On pose
$\Ind_{P}^{G}\sigma=\ind_P^{G}(\sigma\delta_{P}^{1/2})$, où
$\delta_{P}$ est le caractère module de $P$, c'est-à-dire, le
caractère de $M=P/N$ donné par: $\delta_P(m)=[mN_0m^{-1}:N_0]$
pour un arbitraire sous-groupe ouvert compact  $N_0$ de $N$.

(b) Soit $(\pi,V)\in\rRep G$, on note $V_{N}$ le quotient de $V$ par
le sous-espace $V(N)$ engendré par les éléments $\pi(n)x-x$
($n\in N$, $x\in V$). On définit $(r_{P}^{G}(\pi),V_{N})\in \rRep
M$ par
\[r_P^G(m)(v+V(N))=\delta_{P}^{-1/2}(m)(\pi(m)v+V(N)),\ \ m\in M,\ v\in V. \]

\subsection{Rappels de quelques constructions}
Rappelons que dans \cite{BS}, $\S$4, sont  définies deux
catégories $\MOD_{L'/L}$ et $\WD_{L'/L}$:

(i) $\WD_{L'/L}$:  la catégorie des $K$-représentations $(r,N,V)$
du groupe de Weil-Deligne de $L$ (\cite{De}, $\S$8) sur un
$K$-espace vectoriel $V$ de dimension finie telles que la
restriction de $r$ à $W(\overline{\Q}_p/L')$ est non-ramifiée;

(ii) $\MOD_{L'/L}$ (ou $\MOD_{L'/L,K}$ s'il y a risque de confusion
sur $K$): la catégorie des quadruples $(\varphi,N,\Gal(L'/L),D)$
constitués par:
\begin{itemize}
\item[--]  un $L_0'\otimes_{\Q_p}K$-module D libre de rang fini;

\item[--] un Frobenius $\varphi:D\ra D$, c'est-à-dire  une
bijection  $\varphi$ tel que $\varphi((l\otimes k)\cdot
d)=(\varphi'_0(l)\otimes k)\cdot\varphi(d)$, pour $l\in L_0'$,
$k\in K$, $d\in D$;

\item[--]  un endomorphisme $L_0'\otimes_{\Q_p}K$-linéaire $N:D\ra D$
tel que $N\varphi=p\varphi N$;

\item[--] une action de $\Gal(L'/L)$ sur $D$ commutant avec celles de
$\varphi$ et $N$, telle que $g((l\otimes k)\cdot d)=(g(l)\otimes
k)\cdot g(d)$, pour  $g\in \Gal(L'/L)$, $l\in L_0'$, $k\in K$, $d\in
D$.
\end{itemize}
Rappelons aussi que dans \cite{BS}, $\S$4 (ou \cite{Fo2}), sont
définis un foncteur
\[\WD:\MOD_{L'/L}\ra\WD_{L'/L}\]
et un quasi-inverse $\MOD$ de $\WD$. Ces foncteurs induisent une
équivalence entre les deux catégories. On rappelle la construction
dans la suite.
\begin{itemize}
\item[--] Soit $(\varphi,N,\Gal(L'/L),D)$ un objet de $\MOD_{L'/L}$.
Choisissons un plongement $\sigma_0':L_0'\hookrightarrow K$ et
posons $V=D_{\sigma_0'}$. Alors $N$ induit un endomorphisme
$K$-linéaire nilpotent  sur $V$ que l'on note encore par $N$. Si
$w\in W(\overline{\Q}_p/L)$, on définit
$r(w)=\overline{w}\circ\varphi^{-\alpha(w)}$ où\ $\overline{w}$
désigne l'image de $w$ dans $\Gal(L'/L)$ et $\alpha(w)\in f\Z$ est
l'unique entier  tel que l'action induite de $w$ sur $\bFp$ soit la
$\alpha(w)$-puissance du Frobenius arithmétique $x\mapsto x^p$. On
vérifie que $r(w)$ est $L_0'\otimes_{\Q_p}K$-linéaire, et donc
induit un morphisme $K$-linéaire $r(w):V\ra V$. Ceci définit un
objet de $\WD_{L'/L}$. Remarquons que la représentation $(r,N,V)$
est indépendante du choix de $\sigma_0'$ mais à isomorphisme non
canonique près (cf. \cite{BM}, lemme 2.2.1.2).

\item[--] Soit $(r,N,V)$ un objet de $\WD_{L'/L}$ et choisissons un plongement
$\sigma_0':L_0'\hookrightarrow K$.  On pose
\[D=\bigoplus_{n=0}^{f'-1}V_{\sigma_0'\circ\varphi_0'^{-n}},\]
 où
$V_{\sigma_0'\circ\varphi_0'^{-n}}=V$ est muni de l'action de $L_0'$
via $\sigma_0'\circ \varphi_0'^{-n}$, ce qui fait de $D$  un
$L_0'\otimes_{\Q_p}K$-module libre. On définit $\varphi$, $N$ par
\[\left\{ {\begin{array}{ll}
\varphi|_{V_{\sigma_0'\circ\varphi_0'^{1-f'}}}
=r(\omega): V_{\sigma_0'\circ\varphi_0'^{1-f'}}\ra V_{\sigma_0'} &\\

\varphi|_{V_{\sigma_0'\circ\varphi_0'^{-n}}}
 =\id: V_{\sigma_0'\circ\varphi_0'^{-n}}\ra
 V_{\sigma_0'\circ\varphi_0'^{-n-1}}& \mathrm{si\ }0\leq n\leq f'-2 \\\end{array}}\right.\]
où $\omega$ est n'importe quel Frobenius géométrique dans
$W(\overline{\Q}_p/L')$, et
\[\left\{ {\begin{array}{ll}
 N|_{V_{\sigma_0'}}=N|_V:V_{\sigma_0'}\ra V_{\sigma_0'}&  \\

N|_{V_{\sigma_0'\circ\varphi_0'^{-n}}}=p^n\varphi^n\circ
N|_V\circ\varphi^{-n}:V_{\sigma_0'\circ\varphi_0'^{-n}}
\ra V_{\sigma_0'\circ\varphi_0'^{-n}} &\mathrm{si\ }1\leq n\leq f'-1.\\
\end{array}}\right.\]
Finalement, pour $g\in \Gal(L'/L)$, soit $w\in W(\overline{\Q}_p/L)$
un relèvement de $g$, on définit l'action de $g$ sur
$V_{\sigma_0'\circ\varphi_0'^{-n}}$ par
\[g=r(w)\circ\varphi^{\alpha(w)}:V_{\sigma_0'\circ\varphi_0'^{-n}}
\ra V_{\sigma_0'\circ\varphi_0'^{-n-\alpha(w)}},\] où si
$n+\alpha(w)>f'-1$ ou $n+\alpha(w)<0$, on regarde
$V_{\sigma_0'\circ\varphi_0'^{-n-\alpha(w)}}$ comme
$V_{\sigma_0'\circ\varphi_0'^{-k}}$ où $k$ est l'unique entier tel
que
\[0\leq k\leq f'-1 \ \textrm{ et }\ \ n+\alpha(w)\equiv k(\bmod
f').\]

On vérifie que $(\varphi,N,\Gal(L'/L),D)$ ainsi défini est un
objet de $\MOD_{L'/L}$.
\end{itemize}
\begin{rem}\label{rem-scalaire}
Soit $(\varphi,N,\Gal(L'/L),D)$ un objet de $\MOD_{L'/L}$ qui est
absolument irréductible (donc $N=0$). Comme $\varphi^{f'}$ est
$L_0'\otimes_{\Q_p}K$-linéaire et commute avec $\Gal(L'/L)$ et
$\varphi$, on voit que $\varphi^{f'}$ est scalaire à valeur dans
$K^{\times}$.
\end{rem}
\

Maintenant, soit  $(\varphi,N,\Gal(L'/L),D)$ un objet de
$\MOD_{L'/L}$. On définit
\[t_N(D)=\frac{1}{[L:L_0]f'}\val_L(\det{}_{L_0'}(\varphi^{f'}|_D)).\]
Posons $D_{L'}=D\otimes_{L_0'}L'$ et pour $\sigma:L\hookrightarrow
K$,
\[D_{L',\sigma}=D_{L'}\otimes_{L'\otimes_{\Q_p}K}(L'\otimes_{L,\sigma}K).\]
On a alors $D_{L'}\simeq \prod_{\sigma:L\hookrightarrow
K}D_{L',\sigma}$. Donc la donnée d'une filtration décroissante
 exhaustive séparée par des sous-$L'\otimes_{\Q_p}K$-modules
$(\Fil^iD_{L'})_{i}$ (pas forcément libres) stables sous l'action
de $\Gal(L'/L)$ équivaut à la donnée, pour tout $i\in \Z$ et tout
$\sigma:L\ra K$, d'un sous-$L'\otimes_{L,\sigma}K$-module
$\Fil^iD_{L',\sigma}$
 de $D_{L',\sigma}$, qui est stable par
$\Gal(L'/L)$ (donc nécessairement libre), et vérifie
$\Fil^{i+1}D_{L',\sigma}\subset \Fil^iD_{L',\sigma}$ pour tout $i$
et $\sigma$, et $\cup_{i\in\Z}\Fil^iD_{L',\sigma}=D_{L',\sigma}$,
$\cap_{i\in\Z}\Fil^iD_{L',\sigma}=0$ pour tout $\sigma$. Soit
$(\Fil^i D_{L',\sigma})_{i,\sigma}$  une telle filtration. On
définit
\[t_{H}(D_{L'})=\summ_{i\in\Z}\summ_{\sigma}i\dim_{L'}(\Fil^iD_{L',\sigma}/\Fil^{i+1}D_{L',\sigma}).\]
 La filtration est dite \emph{admissible}
si $t_H(D_{L'})=t_N(D)$ et si pour tout sous-$L_0'$-espace vectoriel
$D'$ de $D$ stable par $\varphi$ et $N$, on a $t_H(D'_{L'})\leq
t_N(D')$, où $D'_{L'}$ est muni de la filtration induite par
$D_{L'}$. Par \cite{BM}, proposition 3.1.1.5 et \cite{Fo1},
proposition 4.4.9, pour que la filtration soit admissible, il suffit
de vérifier $t_H(D_{L'})=t_N(D)$ et l'inégalité $t_H(D'_{L'})\leq
t_N(D')$ pour les $L_0'\otimes_{\Q_p}K$-modules $D'$ stables par
$\varphi$, $N$ et $\Gal(L'/L)$ (nécessairement libre),
c'est-à-dire les
sous-objets de $D$ dans $\MOD_{L'/L}$.\\

Si $n\in\N$, on pose $G_n=GL_n(L)$ et $G=G_{d+1}$ où $d\geq 1$ est
un entier fixé. Une \emph{partition} de l'entier $n$ est une suite
$\alpha=(n_1,...,n_r)$ d'entiers positifs tels que
$n=n_1+\cdots+n_r$. Étant donné une telle partition, on note
$P_{\alpha}$ le groupe formé des matrices inversibles triangulaires
supérieures par blocs, c'est-à-dire,
\[P_{\alpha}=\begin{pmatrix}{G_{n_1}}&{*}&{*}&{*}\\{}&{G_{n_2}}&{*}&{*}
\\{}&{}&{\ddots}&{\cdots}\\{ }&{ }&{ }&{G_{n_r}}\end{pmatrix}.\]
On note $N_{\alpha}=N_{P_{\alpha}}$ son radical unipotent et
$M_{\alpha}=M_{P_{\alpha}}$ son quotient de Levi. Les groupes
$P_{\alpha}$ seront appelés
sous-groupe \emph{paraboliques standards} de $G_n$. %Par blocs de
%$\alpha$ on dé signe les ensembles $\{I_i\}$
%\[I_1=\{1,\dots,n_1\},\dots,I_r=\{n_1+\cdots+n_{r-1}+1,n\}.\]
Si $\alpha$ et $\beta$ sont deux partitions de $n$, on dit que
$\beta$ est une sous-partition de $\alpha$ si $M_{\beta}\subseteq
M_{\alpha}$ (notation: $\beta\leq \alpha$).
 %Si $\alpha$ et $\beta$
%sont deux partitions de $n$ de blocs $\{I_i\}$ et $\{J_j\}$

On fixe un choix de  $q^{1/2}$ dans $\overline{\Q}_p$. \`{A} une
représentation $(r,N,V)$ de dimension $d+1$ telle que $r$ soit
semi-simple, la correspondance de Langlands locale permet
d'associer une représentation irréductible lisse $\pi^{\unit}$ de
$G$ sur $\overline{\Q}_p$ normalisée de sorte que le caractère
central de $\pi^{\unit}$ soit $\det(r,N,V)\circ\rec^{-1}$. Notons
que la définition dépend du choix de $q^{1/2}$.  Dans \cite{BS},
cette correspondance est modifiée comme suit:

\begin{itemize}
\item[--] si $\pi^{\unit}$ est générique (cf. \cite{Ku}, \S2.3), alors
$\pi^{\unit}\otimes_{\overline{\Q}_p}|\det|_L^{-d/2}$ admet un
unique modèle sur $K$ qui ne dépend pas du choix de $q^{1/2}$; on
le note $\pi$;

\item[--] si $\pi^{\unit}$ n'est pas générique,  elle s'écrit comme l'unique quotient d'une induction parabolique:
\[\Ind_{P_{\alpha}}^GL(b_1,\tau_1)\otimes\cdots\otimes L(b_s,\tau_s)\]
où $\alpha=(b_1n_1,\dots,b_sn_s)$ est une partition de $d+1$, les
$\tau_i$ sont des représentations irréductibles cuspidales de
$\GL_{n_i}(L)$, les $L(b_i,\tau_i)$ sont des Steinberg
généralisées, c'est-à-dire, $L(b_i,\tau_i)$ est l'unique
quotient irréductible de
\[\Ind_{P_{\alpha_i}}^{G_{b_in_i}}\tau_i\otimes\tau_i|\det|_L\otimes\cdots\otimes\tau_i|\det |_L^{b_i-1}\]
où $\alpha_i$ est la partition $(n_i,\dots,n_i)$ de $b_in_i$
%\[Q=\GL_{d_1}(L)\times\cdots\times\GL_{b_sn_s}(L)\]
 (cf. \cite{Cl}, \S
3.1). On a le résultat suivant:
\\
\begin{lemma} (\cite{BS}, Lemma 4.2)
La représentation
\[(\Ind_{P_{\alpha}}^{G}L(b_1,\tau_1)\otimes\cdots\otimes L(b_s,\tau_s))\otimes_{\overline{\Q}_p}|\det|_L^{-d/2}\]
admet un  modèle unique sur $K$ qui ne dépend pas du choix de
$q^{1/2}$. On le note $\pi$.
\end{lemma}
\end{itemize}
\bigskip
 Si on suppose que $K$ est gros au sens où tous les
$L(b_i,\tau_i)$ sont fixées par $\Gal(\overline{\Q}_p/K)$, alors
la définition de $\pi$ est plus directe. En fait, si on pose
\[\mathcal{L}(b_i,\tau_i)=L(b_i,\tau_i)\otimes_{\overline{\Q}_p}|\det|_L^{(1-b_in_i)/2},\]
alors cette $\mathcal{L}(b_i,\tau_i)$ ainsi définie ne dépend pas
du choix de $q^{1/2}$, et est fixée par $\Gal(\overline{\Q}_p/K)$.
D'après \cite{Cl}, proposition 3.2, elle admet un modèle unique
qui est défini sur $K$ et on le note $\pi_i$. Comme on peut
récrire la représentation originale sous la forme
\[\ind_{P_{\alpha}}^G\mathcal{L}(b_1,\tau_1)\otimes\mathcal{L}(b_2,\tau_2)|\det|_L^{-b_1n_1}
\otimes\cdots\otimes\mathcal{L}(b_s,\tau_s)|\det|_L^{-\sum_{j=1}^{s-1}b_jn_j},\]
 on obtient le modèle $\pi$ en posant
\[\pi=\ind_{P_{\alpha}}^G\pi_1\otimes\pi_2|\det|_L^{-b_1n_1}\otimes\cdots\otimes\pi_s|\det|_L^{-\sum_{j=1}^{s-1}b_jn_j}.\]
%Notons $\chi_i$ le caractè re central de $\pi_i$, alors celui de
%$\pi$ est
%\[\chi=\chi_1\otimes\chi_2|\det|_L^{-b_1n_1}\otimes\cdots\otimes\chi_s|\det|_L^{-\sum_{j=1}^{s-1}b_jn_j}.\]
\subsection{Rappels sur la conjecture et compléments }

%Supposons que $K$ est suffisamment gros pour pouvoir écrir
%\[(r,N,V)=\bigoplus_{i=1}^s(r_i,N_i,V_i)\]
%où\ les $(r_i,N_i,V_i)$ sont toutes  absolument indé
%composables. Soient $(\varphi_i,N_i,\Gal(L'/L),D_i)$ les objets de
%$\MOD_{L'/L}$ correspondant aux $(r_i,N_i,V_i)$. On a $D_i$ est
%absolument indécomposable dans $\MOD_{L'/L}$.

On fixe:

(i) un objet $(r,N,V)$ de $\WD_{L'/L}$ de dimension $d+1$ tel que
$r$ soit semi-simple;

(ii) pour tout $\sigma:L\hookrightarrow K$, un ensemble de $d+1$
entiers $i_{1,\sigma}<\cdots<i_{d+1,\sigma}$.

\`{A}  $(r,N,V)$ comme dans (i) on associe une représentation lisse
$\pi$ comme au \S2.2. Pour des $i_{j,\sigma}$ comme dans (ii), on
pose
\[a_{j,\sigma}=-i_{d+2-j,\sigma}-(j-1),\]
et on note $\rho$ l'unique représentation $\Q_p$-rationnelle de $G$
dont le plus haut poids est $\psi:\diag(x_1,\dots,x_{d+1})\mapsto
\prod_{j=1}^{d+1}\prod_{\sigma:L\hookrightarrow
K}x_j^{a_{j,\sigma}}$ vis-à-vis du sous-groupe  des matrices
triangulaires inférieures (cf. \cite{BS}, \S2).

Si $(r,N,V)\in \WD_{L'/L}$, on note $(r,N,V)^{\textrm{ss}}\in
\WD_{L'/L}$ sa $F$-semisimplification (\cite{De}, $\S$8.5).

%$\mathbf{G}_{\sigma}$ avec poids de plus haut soit
%$(a_{1,\sigma},\cdots,a_{d+1,\sigma})$ vis-à-vis le sous-groupe
%parabolique des matrices triangulaires inférieures, et $\psi$ le
%poids plus haut de $\rho$. Posons $\rho$ l

La conjecture ci-dessous est proposée dans \cite{BS}, Conjecture
4.3:
\begin{conj}\label{conjecture}
Les deux conditions suivantes sont équivalentes:

(i) $\rho\otimes\pi$ admet une norme invariante;

(ii) il existe un objet $(\varphi,N,\Gal(L'/L),D)$ dans
$\MOD_{L'/L}$ tel que:
\[WD(\varphi,N,\Gal(L'/L),D)^{\mathrm{ss}}=(r,N,V),\]
et une filtration admissible $(\Fil^i D_{L',\sigma})_{i,\sigma}$
stable par $\Gal(L'/L)$ sur $D_{L'}$ telle que
\[\Fil^i D_{L',\sigma}/\Fil^{i+1}D_{L',\sigma}\neq 0\Leftrightarrow i\in\{i_{1,\sigma},\dots,i_{d+1,\sigma}\}.\]
\end{conj}

\begin{prop}
La condition (ii) ne dépend pas du choix de $L'$. Plus
précisément, si $L''$ est une autre extension finie de $L$ telle
que
\[ \ [L_0'':\Q_p]=|\Hom(L_0'',K)|\ \ \ \textrm{et} \ \  \ (r,N,V)\in\WD_{L''/L},\]
alors,  (ii) est vraie pour $L'$ si et seulement si (ii) est vraie
pour $L''$.
\end{prop}
\begin{proof}
On peut supposer que $L'\subset L''$.

(a) Soit $(\varphi',N',\Gal(L'/L),D')$ un objet de $\MOD_{L'/L}$. On
pose
\[D''=L_0''\otimes_{L_0'}D',\]
et définit $\varphi''$, $N''$ et l'action de $\Gal(L''/L)$ naturellement. %comme
%suivante:
%\begin{itemize}
%\item[---] $\varphi'':D''\ra D''$ est donnè\ par $\varphi''(l''\otimes
%d')=\varphi''_0(l'')\otimes\varphi'(d')$;
%
%\item[---] $N''=1\otimes N':D''\ra D''$ qui est $L_0''\otimes_{\Q_p}K$-linéaire;
%
%\item[---] pour $g\in \Gal(L''/L)$, on dé finit
%$g(l''\otimes d')=g(l'')\otimes\bar{g}(d')$, où\  $\bar{g}$ est
%l'image de $g$ dans $\Gal(L'/L)$, et $l''\in L''$, $d'\in D'$.
%\end{itemize}
Alors, on vérifie facilement que
%\begin{itemize}
%\item[---] $N''\varphi''=p\varphi''N''$;
%
%\item[---] l'action de $\Gal(L''/L)$ commute avec celui de $\varphi''$
%et $N''$.
%\end{itemize}
%Donc
$(\varphi'',N'',\Gal(L''/L),D'')$ est un objet de $\MOD_{L''/L}$.

De plus, on suppose que
\[\mathrm{WD}(\varphi',N',\Gal(L'/L),D')^{\mathrm{ss}}=(r,N,V)\]
et qu'il existe une filtration $(\Fil^iD'_{L',\sigma})_{i,\sigma}$
sur $D'_{L'}$ comme dans (ii). On munit $D''$ de la filtration
induite:
$\Fil^iD''_{L'',\sigma}=L''\otimes_{L'}\Fil^iD'_{L',\sigma}$. On
vérifie que cette filtration est stable par $\Gal(L''/L)$ et est
admissible, en fait on a $t_N(D')=t_N(D'')$ et
$t_H(D'_{L'})=t_H(D''_{L''})$.
%on a
%\begin{itemize}
%\item[--] cette filtration est pré servé e par $\Gal(L''/L)$;
%
%\item[--] $\Fil^iD''_{L'',\sigma}/\Fil^{i+1}D''_{L'',\sigma}\neq0 \Leftrightarrow
%i\in\{i_{1,\sigma},\cdots,i_{d+1,\sigma}\}$, et donc on déduit
%que $t_H(D'')=t_H(D')$.
%
%\item[--] Par dé finition,
%\[\begin{array}{rll}t_N(D')&=&\frac{1}{[L:L_0]f'}\mathrm{val}_L(\det_{L_0'}(\varphi^{f'}|_{D'}))\\
%&=&\frac{1}{[L:L_0]f''}\mathrm{val}_L(\det_{L_0'}(\varphi^{f''}|_{D'}))\
%\ (\mathrm{comme\ }f'|f'')
%\end{array}\]
%et
%\[t_N(D'')=\frac{1}{[L:L_0]f''}\mathrm{val}_L(\mathrm{det}_{L''_0}(\varphi^{f''}|_{D''})),\]
%donc $t_N(D')=t_N(D'')$.
%\end{itemize}
On fixe un plongement $\sigma_0':L_0'\hookrightarrow K$ et
$\sigma_0'':L_0''\hookrightarrow K$ qui prolonge $\sigma_0'$. Alors
$x\otimes k\mapsto (1\otimes x)\otimes k$ induit un isomorphisme
entre $D'_{\sigma_0'}$ et $D''_{\sigma_0''}$, et cet isomorphisme
commute à $r$ et $N$. Donc on a $\WD(D')\simeq \WD(D'')$.

 (b) Soit $(\varphi'',N'',\Gal(L''/L),D'')\in
\MOD_{L''/L}$ et supposons qu'il existe une filtration
$(\Fil^iD''_{L'',\sigma})_{i,\sigma}$ sur $D''_{L''}$ comme dans
(ii). Posons $D'=D''^{\Gal(L''/L')}$ que l'on munit des structures
induites:  action de $\Gal(L'/L)$, opérateurs $\varphi'$, $N'$ et
filtration.  On vérifie que $(\varphi',N,\Gal(L'/L),D')$ satisfait
à (ii) et
% En fait, ceci est une conséquence de (b) puisque
on a $L_0''\otimes_{L_0'}D''^{\Gal(L''/L')}=D'' $ par Hilbert 90.
\end{proof}

\begin{prop} \label{rem-2-conj}Soit $K'$ une extension galoisienne finie de $K$ contenue dans
$\overline{\Q}_p$. Alors

(1) la condition (i) est vraie pour $K$ si et seulement si (i) est
vraie pour $K'$;

(2) la condition (ii) est vraie pour $K$ si et seulement si (ii) est
vraie pour $K'$ et  la filtration peut être choisie de telle sorte
qu'elle soit stable par l'action de $\Gal(K'/K)$.
\end{prop}
\begin{proof}
(1) est immédiat.

(2) La condition est bien nécessaire. Soit
$(\varphi',N',\Gal(L'/L),D')$ un objet de $\MOD_{L'/L,K'}$
vérifiant (ii) où $D'$ est un $L_0'\otimes_{\Q_p}K'$-module.
Posons $D=D'^{\Gal(K'/K)}$ qui est un $L_0'\otimes_{\Q_p}K$-module.
Comme $\varphi'$, $N'$ et $\Gal(L'/L)$ commutent avec $\Gal(K'/K)$,
on obtient un objet induit $(\varphi,N,\Gal(L'/L),D))$ de
$\MOD_{L'/L,K}$. De plus, puisque la filtration sur $D'_{L'}$  est
stable par $\Gal(K'/K)$, elle induit une filtration sur $D_{L'}$ qui
est bien admissible. Comme
\[\WD(D)^{\textrm{ss}}\otimes_KK'\simeq\WD(D\otimes_KK')^{\textrm{ss}}=(r,N,V)\otimes_KK',\]
on a $\WD(D)^{\textrm{ss}}=(r,N,V)$.
\end{proof}
\begin{rem}
On va voir plus loin que (cf. \S4) si (ii) est vraie pour $K'$,
alors on peut toujours choisir une autre filtration admissible qui
est stable par $\Gal(K'/K)$. Donc en fait (ii) est vraie pour $K'$
si et seulement si (ii) est vraie pour $K$.
\end{rem}

\subsection{L'ordre des $(r_i,N_i,V_i)$}
On conserve les notations du paragraphe précédent. Supposons  $K$
 suffisamment gros pour pouvoir écrire
\[(r,N,V)=\bigoplus_{i=1}^s (r_i,N_i,V_i)\]
avec les $(r_i,N_i,V_i)$ absolument indécomposables de dimension
$d_i$. Soit $(\varphi_i,N_i,\Gal(L'/L),D_i)$ l'objet de
$\MOD_{L'/L}$ tel que $\WD(D_i)=(r_i,N_i,V_i)$, alors $D_i$ est
absolument indécomposable dans $\MOD_{L'/L}$. Posons
$D_{i,0}=\Ker(N_i:D_i\ra D_i)$ et
$\varphi_{i,0}=\varphi_i|_{D_{i,0}}$, alors (notons que $r_i$
correspond au $L(b_i,\tau_i)$, cf. \S2.2)
\[D_i=D_{i,0}\oplus D_{i,0}(1)\oplus \cdots\oplus
D_{i,0}(b_i-1),\] où $D_{i,0}(n)=D_{i,0}$   avec
$\varphi_i|_{D_{i,0}(n)}=p^n\varphi_{i,0}$, et où
$N_i|_{D_{i,0}}=0$ et $N_i$ envoie $D_{i,0}(n)$ dans
$D_{i,0}(n-1)$ par l'identité pour $n>0$. On dit que $D_i$ et
$D_j$ sont \emph{de m\ec me type} si
\[\Hom_{\varphi,\Gal(L'/L)}(D_i,D_j)\neq 0.\] Remarquons que si $D_i$
et $D_j$ sont de m\ec me type, alors il existe un entier $l\in\Z$
tel que $D_{i,0}\simeq D_{j,0}(l)$ comme
$(\varphi,\Gal(L'/L))$-modules.

 On choisit un ordre des $(r_i,N_i,V_i)$ de telle
sorte que (toujours possible):
\[\begin{array}{rll}\{1,\dots,s\}&=&H_1\sqcup H_2\sqcup\cdots\sqcup H_v\\
&=&\{1,\dots,s_1\}\sqcup\{s_1+1,\dots,s_2\}\sqcup\cdots
\sqcup\{s_{r-1}+1,\dots,s\}
\end{array}\] vérifiant les conditions suivantes:
\begin{itemize}
\item[--] $k_1,k_2\in H_i$ si et seulement s'il existe
$j_1,\dots,j_m$ tels que $D_{k_1}$ et $D_{j_1}$, $D_{j_i}$ et
$D_{j_{i+1}}$ ($1\leq i\leq m-1$), $D_{j_m}$ et $D_{k_2}$ sont de
m\ec me type;

\item[--] pour $i$ fixé, $k_1, k_2\in H_i$, $k_1\neq k_2$ et $D_{k_1,0}(l)\simeq D_{k_2,0}$ avec
$l>0$ impliquent $k_1<k_2$;

\item[--]  si $D_{k_1,0}\simeq D_{k_2,0}$, et $b_{k_1}< b_{k_2}$, alors
$k_1<k_2$;
% si $r_{k_1}=r_{k_2}$, on a alors $D_{k_1}\simeq
%D_{k_2}$, on pose $k_1<k_2$ ou $k_1>k_2$;

\item[--] Pour $i\neq j$, posons $M_i=\oplus_{k\in H_{i}}D_k$,
$M_j=\oplus_{k\in H_j}D_k$. Si $t_N(M_i)/\dim M_i< t_N(M_j)/\dim
M_j$, alors $i<j$.
% si $t_N(M_i)/\dim M_i= t_N(M_j)/\dim M_j$, on
%fait n'importe quel choix.
\end{itemize}

\begin{prop}
Choisissons un ordre sur les $(r_i,N_i,V_i)$ comme précédemment,
et soient $L(b_i,\tau_i)$ les représen\-tations définies dans
$\S$2.2, alors la condition ``ne précède pas''  est vérifiée
(\cite{Ku}, Definition 1.2.4).
\end{prop}
\begin{proof}
Supposons l'énoncé faux, et soit $i<j$ tel que $L(b_i,\tau_i)$
précède $L(b_j,\tau_j)$, ce qui implique
$\tau_i(l)=\tau_i|\det|_L^{l}\simeq \tau_j$ pour un entier $l>0$
tel que $l+b_j>b_i$. Alors par la correspondance
\[L(b_i,\tau_i)\leftrightarrow (r_i,N_i,V_i)\leftrightarrow D_i,\]
plus pr\eb cis\eb ment,
\[\tau_i|\det|_L^l\leftrightarrow D_{i,0}(b_i-1-l), \ \ \ \tau_j\leftrightarrow D_{j,0}(b_j-1),\]
on obtient
%, de la relation entre $L(b_i,\tau_i)$ et
%$L(b_j,\tau_j)$, que
\[D_{i,0}\simeq D_{j,0}((l+b_j-b_j)),\]
ce qui contredit $i<j$ puisque $l+b_j-b_i>0$.
\end{proof}

\subsection{Module de Jacquet et un lemme
d'Emerton}\label{section-Emerton} On rappelle un lemme d'Emerton. On
renvoie le lecteur à \cite{Em1} ou \cite{Em3} pour la définition
d'une représentation localement analytique.

Soit $P$ un sous-groupe parabolique de $G$ avec $N$ son radical
unipotent et $M$ son quotient de Levi. Dans \cite{Em3}, \S3, est
défini un foncteur ``module de Jacquet'' $J_P$ ayant la propriété
suivante:
\begin{prop}\label{prop-Emerton}
Soient $\rho$ une $K$-représentation algébrique irréductible de
$G$ et $\pi$ une $K$-représentation admissible lisse de $G$, alors
il existe un isomorphisme canonique de $M$-représentations
\[J_P(\rho\otimes_K\pi)\simto \rho^N\otimes_K(r_P^G\pi)\delta_P^{1/2}.\]
\end{prop}
\begin{proof}
Voir \cite{Em3}, Proposition 4.3.6.
\end{proof}

Si $N_0$ est un sous-groupe ouvert compact de $N$, on note $Z(M)$ le
centre de $M$ et on pose
\[Z(M)^+=\{z\in Z(M)|\ zN_0z^{-1}\subset N_0\}.\]

\begin{lemma}(Emerton)\label{lemme-Emerton}
Avec les notations ci-dessus. Soient $V$ une représentation
localement analytique de $G$ et $\chi:Z(M)\ra K^{\times}$ un
caractère localement analytique de $Z(M)$. On considère les deux
conditions:
\begin{itemize}
\item[(1)] $V$ admet une norme invariante;

\item[(2)] Si $\Hom_{Z(M)}(\chi,J_P(V))\neq 0$, alors
\[|\chi(z)\delta_P^{-1}(z)|_p\leq 1,\ \ \ \ \forall z\in Z(M)^{+}.\]
\end{itemize}
Alors (1) entraîne (2).
\end{lemma}
\begin{proof}
Voir \cite{Em2}, Lemma 1.6 ou \cite{Em3}, Lemma 4.4.2.
\end{proof}
\begin{defn}\label{def-cond-de-Emerton}
On dit que $V$ satisfait à la \emph{condition d'Emerton} si pour
tout
  $P$,  $N_0$ et
 $\chi$ comme ci-dessus la condition (2) du lemme \ref{lemme-Emerton} est
 vérifiée.
\end{defn}
\begin{rem}\label{rem-cond-de-Emerton}
Comme tout sous-groupe parabolique de $G$ est conjugué à un
parabolique standard, on peut supposer que $P$ est un parabolique
standard dans la définition \ref{def-cond-de-Emerton}. Aussi on
peut supposer $N_0=N\cap \textrm{M}_{d+1}(\Z_p)$, où
$\textrm{M}_{d+1}(\Z_p)$ est l'anneau des matrices carrées
$(d+1)\times (d+1)$ à coefficients dans $\Z_p$.
\end{rem}
\subsection{Théorème principal}

Le résultat principal de cet article est le suivant.
\begin{theorem}\label{thm-main}
Supposons  $K$  suffisamment gros tel que l'on a
\[(r,N,V)=\bigoplus_{i=1}^s(r_i,N_i,V_i)\]
avec les $r_i$ absolument indécomposables de dimension $d_i$. On
considère les quatre conditions suivantes:

(i) $\rho\otimes\pi$ admet une norme invariante;

(ii) Il existe un objet $(\varphi,N,\Gal(L'/L),D)$ dans
$\MOD_{L'/L}$ tel que:
\[WD(\varphi,N,\Gal(L'/L),D)^{\mathrm{ss}}=(r,N,V),\]
et une filtration admissible $(\Fil^i D_{L',\sigma})_{i,\sigma}$
stable par $\Gal(L'/L) \times\Gal(K/\Q_p)$ sur $D_{L'}$ telle que
\[\Fil^i D_{L',\sigma}/\Fil^{i+1}D_{L',\sigma}\neq 0\Leftrightarrow i\in\{i_{1,\sigma},\dots,i_{d+1,\sigma}\}.\]

(iii) Avec l'ordre des $D_i$ comme en \S2.4, les inégalités
suivantes sont verifiées
%pour tout permutation $\nu\in S_n$:
\[[K:L]\summ_{j=1}^{d_1}\summ_{\sigma}i_{j,\sigma}\leq t_N(D_{1}),\]
\[\vdots\]
\[[K:L]\summ_{j=1}^{d_1+\cdots+d_{s-1}}\summ_{\sigma}i_{j,\sigma}\leq\summ_{i=1}^{s-1}t_N(D_{i}),\]
\[[K:L]\summ_{j=1}^{d+1}\summ_{\sigma}i_{j,\sigma}=\summ_{i=1}^{s}t_N(D_{i})=t_N(D).\]

%(iii)'  quand $\nu=\id$, les iné galité s dans (iii) sont
%verifié es.
%\[[K:L]\summ_{j=1}^{d_1}\summ_{\sigma}i_{j,\sigma}\leq t_N(D_1),\]
%\[\cdots\cdots\]
%\[[K:L]\summ_{j=1}^{d_1+\cdots+d_{r-1}}\summ_{\sigma}i_{j,\sigma}\leq\summ_{i=1}^{r-1}t_N(D_i),\]
%\[[K:L]\summ_{j=1}^{d+1}\summ_{\sigma}i_{j,\sigma}=\summ_{i=1}^{r}t_N(D_i)=t_N(D).\]

(iv)  Le caractère central de $\rho\otimes\pi$ est unitaire  et $V$
satisfait la condition d'Emerton.

 Alors, $(i)\Rightarrow(ii)\Leftrightarrow(iii)\Leftrightarrow(iv)$.
\end{theorem}
\begin{rem}\label{rem-1-main}
(1) L'implication ``$(ii)\Rightarrow(iii)$'' découle de la
définition d'admissibilité.

(2) L'implication ``$(i)\Rightarrow(iv)$'' est une redite du lemme
d'Emerton  et \cite{BS}, Proposition 5.1.

\end{rem}
\begin{rem}\label{rem-2-main}
Si on pose $b=b_1+\cdots+b_s$,
\[D_1'=D_{1,0},\dots,D_{b_1}'=D_{1,0}(b_1-1), \dots,D_b'=D_{s,0}(b_s-1)\]
et $d_i'=\textrm{rg}_{L_0'\otimes_{\Q_p}K}D_i'$, alors pour toute
permutation $\nu\in S_b$, on a
\[[K:L]\summ_{j=1}^{d'_{\nu(1)}}\summ_{\sigma}i_{j,\sigma}\leq t_N(D'_{\nu(1)}),\]
\[\vdots\]
\[[K:L]\summ_{j=1}^{d'_{\nu(1)}+\cdots+d'_{\nu(b-1)}}\summ_{\sigma}i_{j,\sigma}\leq\summ_{i=1}^{b-1}t_N(D'_{\nu(i)}),\]
\[[K:L]\summ_{j=1}^{d+1}\summ_{\sigma}i_{j,\sigma}=\summ_{i=1}^{b}t_N(D'_{\nu(i)})=t_N(D).\]
\end{rem}
\begin{proof}
Ceci résulte de l'ordre des $(r_i,N_i,V_i)$ (cf. \S2.4).
\end{proof}
\begin{rem}
D'après le théorème, la conjecture \ref{conjecture} se
``réduit'' à (iv)$\Rightarrow$(i) (supposons $K$ suffisamment
gros). Il y a un résultat analogue dans le cas complexe (cf.
\cite{Cas}, Theorem 4.4.6).
\end{rem}

\begin{cor}\label{cor-main}
Avec les notations comme dans la conjecture \ref{conjecture}, on a
l'implication
\[(i)\Rightarrow (ii).\]
\end{cor}
\begin{proof}
Ceci est une conséquence du théorème \ref{thm-main} et la
proposition \ref{rem-2-conj}.
\end{proof}
\

Dans toute la suite de cet article, on suppose que $K$ est
suffisamment gros au sens du théorème \ref{thm-main}.

 \section{Preuve de (iii)$\Leftrightarrow$ (iv)}
Dans ce chapitre, on prouve que (iii) entraîne (iv) (sous
l'hypothèse que $K$ est suffisamment gros). Conservons les
notations précédentes: $G=GL_{d+1}(L)$, $d_i=b_in_i$,
$\alpha=(d_1,\dots,d_s)$, $\tau_i$ est une représentation cuspidale
irréductible de $G_{d_i}$, $P=P_{\alpha}$ est le sous-groupe
parabolique standard correspondant à la partition $\alpha$,
$N_{\alpha}$ son radical unipotent, $M_{\alpha}$ son quotient de
Levi, $\pi$ est l'unique modèle sur $K$ de la représentation
\[\Ind_{P_{\alpha}}^GL(b_1,\tau_1)\otimes\cdots\otimes L(b_s,\tau_s)\otimes_{\bQp}|\det|^{-d/2}_L,\]
 $\rho$ est la représentation rationnelle de $G$ associées aux
entiers $\{a_{j,\sigma}\}_{j,\sigma}$ et $\psi$ son plus haut poids.

 \subsection{Le foncteur $r_{Q}^G\circ \Ind_{P}^G$}

  On note %$M=G_{\beta}$ et $N=U_{\beta}$ et
  $\tau=L(b_1,\tau_1)\otimes\cdots\otimes L(b_s,\tau_s)$.  On fixe $ Q=P_{\beta}$
  un parabolique standard de $G$ où $\beta=(m_1,\dots,m_r)$.
Dans ce n$^{\circ}$, on trouve une condition nécessaire pour qu'un
caractère $\chi$ de $Z({M_{\beta}})$ soit tel que
\[\Hom_{Z(M_{\beta})}(\chi,r_{Q}^G\circ\Ind_{P }^G(\tau))\neq 0.\]
%Notons qu'on peut supposer $Q=P_{\beta}$ est un parabolique

D'abord on a besoin du lemme suivant:
\begin{lemma}\label{lemme-Zel}
Soient $l=km$ et $\tau'$ une représentation cuspidale irréductible
de $G_m$. Si $\gamma=(l_1,\dots,l_n)$ est une partition de $l$,
alors $r_{P_{\gamma}}^{G_l}(L(k,\tau'))\neq 0$ si et seulement si
$m|l_i$ pour tout $i$, et dans ce cas
$r_{P_{\gamma}}^{G_l}(L(k,\tau'))$ est irréductible  et égal à
 \[L(p_n,\tau'|\det |^{p_1+\cdots+p_{n-1}})\otimes\cdots\otimes L(p_1,\tau'),\]
 où $p_i=l_i/m$.  Son
caractère central est égal  à la restriction de
\[\big(\chi(\tau')|\det|_L^{p_1+\cdots+p_{n-1}}\otimes\chi(\tau')|\det|_L^{k-1}\big)\otimes
\cdots\otimes\big(\chi(\tau')\otimes\cdots\otimes\chi(\tau')|\det|_L^{p_1-1}\big)
\] à $Z({M_{\gamma}})$, où\ $\chi(\tau')$ désigne le caractère central de $\tau'$.
\end{lemma}
\begin{proof}
C'est une conséquence de la proposition 9.5 et la proposition 1.5
de \cite{Zel}.
\end{proof}

Posons $\gamma=(\underset{b_1\
fois}{\underbrace{n_1,\dots,n_1}},\underset{b_2\
fois}{\underbrace{n_2,\dots,n_2}},\dots,\underset{b_s\
fois}{\underbrace{n_s,\dots,n_s}})$ la partition de $d+1$, et
\[\chi_{\gamma}=(\chi(\tau_1)\otimes\cdots\otimes\chi(\tau_1)|\det|_L^{b_1-1})\otimes\cdots\otimes
(\chi(\tau_s)\otimes\cdots\otimes\chi(\tau_s)|\det|_L^{b_s-1})\]
alors $\chi_{\gamma}$ est un caractère lisse de $Z({M_{\gamma}})$
dont le restriction à $Z(M_{\alpha})$ est le caractère central de
$\tau$.
 Posons $W\simeq S_{d+1}$ le
groupe de Weyl de $G$, $P_0=P_{(1,\dots,1)}$ et
\[W^{\alpha,\beta}=\{w\in W|w(M_{\alpha}\cap P_0)w^{-1}\subset P_0, w^{-1}(M_{\beta}\cap P_0)w\subset P_0\}.\]
Si $w\in W^{\alpha,\beta}$ et on pose $\alpha'=\alpha\cap
w^{-1}\beta w$ (cf. \cite{Zel}, \S1.2), alors $\alpha'\leq\alpha$.
Si de plus $\gamma\leq\alpha'$, c'est-à-dire, $M_{\gamma}\subseteq
M_{\alpha'}$, on a $Z({M_{\alpha'}})\subseteq Z({M_{\gamma}})$,
via lequel on peut voir $\chi_{\gamma}$ comme un caractère de
$M_{\alpha'}$ et de même on voit $\chi_{\gamma}^{w^{-1}}$ comme un
caractère de $M_{\beta'}$, où par définition $\beta'=w\alpha
w^{-1}\cap\beta$ et $\chi^{w^{-1}}(z)=\chi(w^{-1}zw)$.

\begin{prop}\label{prop-Zel}
(1) Si $P=Q$, on a
$\Hom_{Z({M_{\beta}})}(\chi_{\gamma},r_{P}^G\circ\Ind_{P}^G(\tau))\neq
0$.

(2)
$\Hom_{Z({M_{\beta}})}(\chi,r_{Q}^{G}\circ\Ind_{P}^{G}(\tau)\neq0$
si et seulement s'il existe une permutation $w\in W^{\alpha,\beta}$
telle que
\[ \gamma\leq \alpha'=\alpha\cap w^{-1}\beta w\ \textrm{ et }\ \chi=\chi_{\gamma}^{w^{-1}}|_{Z({M_{\beta}})}.\]
\end{prop}
 %et $Q$  le sous-groupe parabolique triangulaire  supé rieur de
%type $(d_1,\cdots,d_s)$ où\ $d_i=d_i$, c'est-à-dire,
%\[Q=\begin{pmatrix}{\GL_{d_1}(L)}&{*}&{*}&{*}\\{}&{\GL_{d_2}(L)}&{*}&{*}
%\\{}&{}&{\ddots}&{\cdots}\\{ }&{ }&{ }&{\GL_{d_s}(L)}\end{pmatrix},\]
 % $N$ son radical nilpotent et $M$ son quotient de
%Levi. Posons $Q^*=w_0Qw_0^{-1}$, $N^*=w_0Nw_0^{-1}$,
%$M^*=w_0Mw_0^{-1}$.

%Si $\beta\leq \alpha$, on dé finit
%$\ind_{G_{\beta}}^{G_{\alpha}}$, $\Ind_{G_{\beta}}^{G_{\alpha}}$,
%$J_{G_{\beta}}^{G_{\alpha}}$, $r_{G_{\beta}}^{G_{\alpha}}$
%naturellement.
\begin{proof}
(1) est un cas spécial de (2), et  (2) se déduit de \cite{Zel},
Theorem 1.2 et Proposition 1.6, et du lemme \ref{lemme-Zel}.
 \end{proof}

Posons \[w_0=\begin{pmatrix}{}&{}&1\\{}&{\iddots}&{}\\
1&{}&{}\end{pmatrix}\in\GL_{d+1}(L) , \] et $Q^*=w_0Qw_0^{-1}$,
$N_{\beta}^*=w_0N_{\beta}w_0^{-1}$,
$M_{\beta}^*=w_0M_{\beta}w_0^{-1}$.

Par la proposition \ref{prop-Zel}, on obtient
\begin{cor}\label{cor-nonnul}
(1) Si $P=Q$, on a
$\Hom_{Z(M_{\beta}^{*})}\big((\chi_{\gamma})^{w_0}\psi\delta_{P^*}^{1/2}|\det|_L^{-d/2},
J_{P^{*}}(\rho\otimes\pi)\big)\neq 0$;

(2)
$\Hom_{Z(M_{\beta}^{*})}\big(\chi,J_{Q^{*}}(\rho\otimes\pi)\big)\neq
0$ si et seulement s'il existe $w\in W^{\alpha,\beta}$ tel que
\[ \gamma\leq
\alpha'=\alpha\cap w^{-1}\beta w,\ \textrm{ et }\
\chi=(\chi_{\gamma}^{w^{-1}}|_{Z(M_{\beta}^*)})^{w_0}\psi\delta_{Q^*}^{1/2}|\det|_L^{-d/2}.\]
\end{cor}
\begin{proof}
Ce corollaire se déduit de la proposition \ref{prop-Zel} et des
faits suivants:

(a) d'après la proposition \ref{prop-Emerton},
\[J_{Q^*}(\rho\otimes\pi)\simeq \rho^{N_{\beta}^*}\otimes
(r_{Q^*}^G\pi)\delta_{Q^*}^{1/2}=\rho^{N_{\beta}^*}\otimes(r_{Q^*}^G\circ\Ind_P^G\tau)|\det|_L^{-d/2}\delta_{Q^*}^{1/2};\]

(b) pour tout $Q$ parabolique standard, $\rho^{N^{*}}$ est une
représentation algébrique irréductible de $M^{*}$ dont le
caractère central est la restriction de $\psi$ à ${Z(M^*})$;

(c) si $\sigma\in \rRep G$, alors
$\Hom_{Z(M_{\beta})}(\chi,r_Q^G(\sigma))\neq0$ si et seulement si
$\Hom_{Z(M_{\beta}^{*})}(\chi^{w_0},\linebreak
r_{Q^{*}}^G(\sigma))\neq 0$.
\end{proof}

 \subsection{La preuve}

\begin{proof}
[$(iii)\Rightarrow (iv)$:] D'abord par \cite{BS}, Proposition 5.1,
le caractère central de $\rho\otimes\pi$ est unitaire.

 Soit
$\beta=(m_1,\dots,m_r)$ une partition de $d+1$. Posons
$Q=P_{\beta}$, $Q^*=w_0Qw_0^{-1}$, $N^*_0=N_{\beta}^*\cap
\mathrm{M}_{d+1}(\Z_p)\subset N_{\beta}^*$,
$Z(M_{\beta}^*)^+=\{z\in Z(M_{\beta}^*)|zN_0^*z^{-1}\subset
N_0^{*}\}$. Soit $w\in W^{\alpha,\beta}$ tel que
$\gamma\leq\alpha'= \alpha\cap w^{-1}\beta w$. D'après le
corollaire \ref{cor-nonnul} et la remarque
\ref{rem-cond-de-Emerton}, il suffit de prouver \[ \val_p\big(
\chi_{\gamma}^{w^{-1}w_0}\psi\delta_{Q^*}^{-1/2}|\det|_L^{-d/2}(z)\big)\geq
0,
 \ \ \forall z\in Z(M_{\beta}^*)^+.  \leqno (*) \]
%Comme $\beta'=\beta\cap w\alpha w^{-1}\leq \beta$ et
%$\chi_{\gamma}^{w^{-1}}$ restreint à\ $G_{\beta'}$,
On récrit  $\chi_{\gamma}^{w^{-1}}$ sous la forme
\[\chi_1'\otimes\cdots\otimes\chi_r'\] avec $\chi'_i: Z(G_{m_i})\ra K^{\times}$.
Soit $t_1\leq t_2\leq \cdots\leq t_r$ et on pose
\[z=z(t_1,\dots,t_r)=\diag(p^{t_1}\ide_{m_r},\dots,p^{t_r}\ide_{m_1}),\]
où $\ide_{m_i}$ désigne la matrice identité d'ordre $m_i$. Alors
$z\in Z(M_{\beta}^*)^+$ et on voit qu'il suffit de prouver (*) pour
$z\in Z(M_{\beta}^*)^+$ de cette forme.

On calcule $\val_p(\chi_{\gamma}^{w^{-1}w_0}(z))$,
$\val_p(\psi(z))$, $\val_p(\delta_{Q^*}^{-1/2}(z))$ et
$|\det|_L^{-d/2}$ respectivement.
\begin{itemize}
 \item[(a)]  Trivialement on a
  \[\val_p(\chi_{\gamma}^{w^{-1}w_0}(z))=\summ_{i=1}^r{t_{s+1-i}}\big(\val_p(\chi_i'(p\ide_{m_i}))
  \big)\]
et
\[\val_p(|\det(z)|_L^{-d/2})=\frac{d}{2}[L:\Q_p]\summ_{i=1}^rm_it_{r+1-i}.\]

\item[(b)] Par la définition de $\rho$, on a
\[\begin{array}{rll}\val_p(\psi(z))&=&t_1\big(\summ_{j=1}^{m_r}\summ_{\sigma}a_{j,\sigma}\big)
+\cdots+t_{r}\big(\summ_{j=m_r+\cdots+m_2+1}^{d+1}\summ_{\sigma}a_{j,\sigma}\big)\\
&=&-t_1\summ_{j=1}^{m_r}\summ_{\sigma}i_{d+2-j,\sigma}-\cdots-t_r\summ_{j=m_r+\cdots+m_2+1}^{d+1}
\summ_{\sigma}i_{d+2-j,\sigma}\\
&&
-[L:\Q_p]t_1\frac{m_r(m_r-1)}{2}-\cdots-t_r[L:\Q_p]\frac{m_1(2\summ_{j>1}m_{j}+m_{1}-1)}{2}\\
&=&-\summ_{i=1}^r(t_{r+1-i}\summ_{j=m_1+\cdots+m_{i-1}+1}^{m_1+\cdots+m_{i}}\summ_{\sigma}i_{j,\sigma})\\
&&-\frac{1}{2}[L:\Q_p]\big(\summ_{i=1}^rt_{r+1-i}m_{i}(2\summ_{j>i}m_{j}+m_{i}-1)\big).\end{array}\]

\item[(c)] Par \cite{Vi1}, chapitre I.2.7, exemple (c), on a
%$\delta_Q(z)=q^{d(z)}$ avec
%\[d(z)=(\summ_{j<i}d_j-\summ_{j>i})p^{t_1}\]
\[\begin{array}{rll}\val_p(\delta_{Q^*}^{-1/2}(z))&=&\frac{1}{2}\val_p(\delta_Q^{-1}(w_0zw_0^{-1}))\\
&=&-\frac{1}{2}[L:\Q_p]\summ_{i=1}^r(\summ_{j<i}m_j-\summ_{j>i}m_j)m_it_{r+1-i}.\end{array}\]

\end{itemize}
Un calcul rapide implique que
\[\begin{array}{rll}&&\val_p(\chi_{\gamma}^{w^{-1}w_0}\psi\delta_{Q^{*}}^{-1/2}\det|_L^{-d/2}(z)|)\\
&=&\summ_{i=1}^rt_{r+1-i}\val_p(\chi_i'(p\ide_{m_i}))-
\summ_{i=1}^r(t_{r+i-1}\summ_{j=m_1+\cdots+m_{i-1}+1}^{m_1+\cdots+m_{i}}\summ_{\sigma}i_{j,\sigma}).\end{array}\]
Notons que, par la correspondance de Langlands locale unitaire
(non modifiée), le caractère central de $L(b_i,\tau_i)$
  est $\det(r_i,N_i,V_i)\circ\rec^{-1}$, on a donc pour tout
$1\leq i\leq s$ (cf. \cite{BS}, Proposition 5.1)
\[%\chi(\tau_i)(p\ide_{m_i})=t_N(D_{i,0}(b_i-1))\ \textrm{ et }\
\chi(\tau_i)|\det|_L^{j}(p\ide_{m_i})=\frac{1}{[K:L]}t_N(D_{i,0}(b_i-1-j)),
\ \ j=0,\dots,b_i-1;\] mais, par d\eb finition, $\chi_i'$ est un
produit de tels caract\ed res, on en déduit les inégalités (*) de
la remarque \ref{rem-2-main} et du fait que $t_1\leq t_2\leq
\cdots\leq t_r$.
%par exemple, en prenant $t_0=\cdots t_{r-1}=0$ et $t_r=1$, on trouve
%\[\val_p(\chi_i'(p\ide_{m_i}))\geq [K:L]\summ_{j=1}^{m_1}\summ_{\sigma}i_{j,\sigma}\].

(iv)$\Rightarrow$(iii): Immédiat à partir du calcul ci-dessus.
\end{proof}
\begin{rem}
Cette preuve est une généralisation de \cite{Em2}, Lemma 2.1.
\end{rem}

\section{Preuve de (iii)$\Rightarrow$(ii)}
Supposons que (iii) est vrai, on  va prouver (ii) dans ce chapitre.
 Posons
\[(\varphi_i,N_i,\Gal(L'/L),D_i)=\MOD(r_i,N_i,V_i).\]

\noindent\emph{Convention:} comme les
$L_0'\otimes_{\Q_p}K$-modules $D'$
 qu'on va traiter sont tous objets de $\MOD_{L'/L}$ libres de rang égal
 à $\dim_K(\WD(D'))$, on pose $\dim D'=\mathrm{rg}_{L_0'\otimes_{\Q_p}K}D'$ .
\subsection{Construction de $(\varphi,N,\Gal(L'/L),D)$}
 Dans ce paragraphe, on construit un  objet
$(\varphi,N,\Gal(L'/L),D)$ de $\MOD_{L'/L}$ tel que
\[\WD(\varphi,N,\Gal(L'/L),D)^{\textrm{ss}}=(r,N,V).\]
D'abord, on regarde quelques exemples.
\begin{examples}\label{examples}
On considère le cas $L'=L=\Q_p$ (donc $\Gal(L'/L)$=\{1\}).
Définissons $(D_0,\varphi_0,N_0)$ par:
\[D_0=Ke,\  \varphi_0 (e)=e,\ N_0(e)=0.\]

(1) (a) Soit $D=D_1\oplus D_2=D_0\oplus (D_0\oplus D_0(1))$,
c'est-à-dire,
\[D=K^3=Ke_1\oplus Ke_2\oplus Ke_3\] avec
\[\varphi(e_1)=e_1,\ \ \varphi (e_2)= e_{2}, \ \ \varphi (e_3)=pe_3\] et \[N(e_1)=0,
\ \ N(e_2)=0, \ \ N(e_3)=e_2.\] Alors il existe un triplet
$(i_1,i_2,i_3)$ tel que
\[i_1<i_2<i_3,\ \ \ i_1+i_2+i_3=t_N(D)=1,\]
et tel qu'il n'existe aucune filtration admissible dont les poids de
Hodge-Tate (i.e. les entiers $i$ tels que
$\Fil^{-i}D/\Fil^{-i+1}D\neq 0$ avec multiplicité) soient
$\{-i_1,-i_2,-i_3\}$. En fait, on peut choisir un triplet
$(i_1,i_2,i_3)$ tel que
\[i_1<i_2<i_3,\ \ \ i_1+i_2+i_3=t_N(D)=1, \textrm{ et }i_2>0.\]
Soit $(\Fil^iD)_i$ une filtration  quelconque dont les poids de
Hodge-Tate
 sont $\{-i_1,-i_2,-i_3\}$. Alors $\dim
\Fil^{i_2}D=2$ et si on pose $D'=\Fil^{i_2}D\cap (Ke_1\oplus
Ke_2)$ qui est un sous-objet de $D$, on a
\[\left\{ {\begin{array}{ll}
t_H(D')=i_2+i_3> 0=t_N(D')
& \textrm{si }\Fil^{i_2}D=Ke_1+Ke_2 \\

t_H(D')\geq i_2>0=t_N(D')& \textrm{sinon}. \\\end{array}}\right.\]

% Mais on a $t_N(D')=1$, donc
%$\Fil$ n'est pas admissible.

La situation sera améliorée si on modifie $\varphi$ par
\[\varphi^{*} (e_1)=e_1+e_2,\ \  \varphi^{*} (e_2)=e_2,\ \ \varphi^{*} (e_3)=pe_3.\]
On vérifie que $p\varphi^* N=N\varphi^*$, et tout sous-objet $D'$
non trivial de $D$ est de la forme:
\[Ke_2;\ \ \ Ke_i\oplus Ke_j,\ (i,j)=(1,2)\textrm{ ou }(2,3).\]
 On peut
construire une filtration telle que \[t_H(Ke_i)=i_1, \ \
t_H(Ke_i\oplus Ke_j)=i_1+i_2\] pour $i\neq j$ (cf. lemme
\ref{lemma-BS} ci-après), et la condition $i_1<i_2<i_3$ entraîne
que cette filtration est admissible.

(b) La même preuve  pour $D=(D_0\oplus D_0(1))\oplus D_0(1)$.

 (2) Soit $D=D_1\oplus D_2=(D_0\oplus D_0(1))\oplus (D_0(1)\oplus
 D_0(2))$, c'est-à-dire,
 \[D=Ke_1\oplus Ke_2\oplus Ke_3\oplus Ke_4,\]
\[\varphi(e_1)=e_1,\ \ \varphi (e_2)=pe_2,\ \ \varphi (e_3)=pe_3, \ \ \varphi (e_4)=p^2e_4;\]
\[N(e_1)=0,\ \ N(e_2)=e_1,\ \ N(e_3)=0,\ \ N(e_4)=e_3.\]
On peut énumérer tous les sous-objets non triviaux de $D$ stables
par $\varphi$, $N$:
\[Ke_i,\ i=1,3;\ \ Ke_i\oplus Ke_j, \  (i,j)=(1,2),(1,3),(3,4);\]
\[Ke_i\oplus Ke_j\oplus Ke_k,\  (i,j,k)=(1,2,3), (1,3,4);\]
\[Ke_1\oplus Kv,\ \ \mathrm{avec\ }   v=ae_2+be_3,\ a,b\neq 0.\]
Soient $\{i_j\}_{1\leq j\leq 4}$ des entiers tels que
\[\summ_{j=1}^4i_j=4,\ \ \ i_1<i_2<i_3<i_4.\]
 On peut construire
une filtration $(\Fil^iD)_i$ sur $D$ telle que (cf. lemme
\ref{lemma-BS})
\[\Fil^iD/\Fil^{i+1}D\neq 0\Leftrightarrow i\in\{i_1,i_2,i_3,i_4\}\] et telle que
\[t_H(Ke_i)=i_1, \ \ t_H(Ke_i\oplus Ke_j)=i_1+i_2,\ \ t_H(Ke_i\oplus Ke_j\oplus Ke_k)=i_1+i_2+i_3,\]
pour tous $i$, $j$, $k$ distincts. En particulier, si
$D'=Ke_1\oplus Kv$ avec $v=ae_2+be_3,\ a,b\neq 0$, on a
$t_H(D')\leq i_1+i_3$. L'hypothèse sur $\{i_j\}_{1\leq j\leq 4}$
implique $i_1+i_3\leq 1$ et on en déduit que la filtration est
admissible.

Notons qu'on peut aussi modifier $\varphi$ comme en (1):
\[\varphi^* (e_2)=pe_2+pe_3, \ \ \varphi^* (e_i)=\varphi e_i \ \mathrm{pour\ }i\neq 2.\]
 On vérifie que $p\varphi^*N=N\varphi^*$ et
\[WD(\varphi,N,\Gal(L'/L),D)^{\mathrm{ss}}=WD(\varphi^*,N,\Gal(L'/L),D)^{\mathrm{ss}},\]
et  que la filtration définie ci-dessus est encore admissible.

(3) Soit $D=D_1\oplus D_2=\big(D_0\oplus D_0(1)\oplus
D_0(2)\big)\oplus D_0(1)$. De fa\c{c}on plus précise,
 \[D=Ke_1\oplus Ke_2\oplus Ke_3\oplus Ke_4;\]
\[\varphi (e_1)=e_1,\ \ \varphi (e_2)=pe_2,\ \ \varphi (e_3)=p^2e_3, \ \ \varphi (e_4)=pe_4;\]
\[N(e_1)=0,\ \ N(e_2)=e_1,\ \ N(e_3)=e_2,\ \ N(e_4)=0.\]
On peut vérifier que si $\varphi^*$ est un endomorphisme de $D$ tel
que $p\varphi^*N=N\varphi^*$ et tel que
\[WD(\varphi^*,N,\Gal(L'/L),D)^{\mathrm{ss}}=WD(\varphi,N,\Gal(L'/L),D)^{\mathrm{ss}},\]
alors $\varphi^*=\varphi$.
\end{examples}
\begin{rem}
(1) L'exemple dans (1) nous dit que l'on \emph{doit} modifier
$\varphi$ dans certains cas. (Notons que le triplet $(i_1,i_2,i_3)$
vérifie la condition (iii) du théorème \ref{thm-main}.)

(2) L'exemple (2) signale qu'il n'est \emph{pas toujours
nécessaire} de modifier $\varphi$.

(3) L'exemple (3) nous dit qu'il n'est \emph{pas toujours possible}
de modifier $\varphi$.

\end{rem}
\begin{lemma}\label{lemma-homo}
On fixe un $H_i$, et soit $k_1\in H_i$, $k_2\in H_i$ tels que $k_1<
k_2$. On peut écrire
\[D_{k_1}=D_0\oplus D_0(1)\oplus\cdots\oplus D_0(r_1),\]
\[D_{k_2}=D_0(l)\oplus D_0(l+1)\oplus\cdots\oplus D_0(l+r_2),\]
avec $D_0=\ker(N:D_{k_1}\ra D_{k_1})$ absolument irréductible et
$l\geq 0$. Alors

(1) Pour que $\Hom_{\MOD_{L'/L}}(D_{k_1},D_{k_2})\neq 0$ il faut et
il suffit que $l\leq r_1\leq l+r_2$. S'il en est ainsi,
$\Hom_{\MOD_{L'/L}}(D_{k_1},D_{k_2}) $ est un $K$-espace vectoriel
de dimension 1.

(2) Supposons $\Hom_{\MOD_{L'/L}}(D_{k_1},D_{k_2})\neq 0$ et soit
$\alpha: D_{k_1}\ra D_{k_2}$ un tel morphisme non nul. On modifie
le Frobenius sur $D_{k_1}\oplus D_{k_2}$ en posant
$\varphi:D_{k_1}\oplus D_{k_2}\ra D_{k_1}\oplus D_{k_2}$,
\[\varphi(e_1)=\varphi_{k_1}(e_1)+\varphi_{k_2}(\alpha(e_1)),\ \ \ \varphi(e_2)=\varphi_{k_2}(e_2),\]
où\ $\varphi_{k_i}$ est le Frobenius sur $D_{k_i}$.  Alors
$p\varphi N=N\varphi$ où\ $N=N_{k_1}\oplus N_{k_2}: D_{k_1}\oplus
D_{k_2}\ra D_{k_1}\oplus D_{k_2}$. De plus, on a
\[(\varphi,N,\Gal(L'/L),D_{k_1}\oplus D_{k_2})\in\MOD_{L'/L},\] et
\[\WD(\varphi,N,\Gal(L'/L),D_{k_1}\oplus D_{k_2})^{\mathrm{ss}}
=\WD(\varphi_{k_1}\oplus \varphi_{k_2},N,\Gal(L'/L),D_{k_1}\oplus
D_{k_2})^{\mathrm{ss}}.\]

(3) Tout sous-objet de $(\varphi,N,\Gal(L'/L),D_{k_1}\oplus
D_{k_2})$ est, \ad\ isomorphisme pr\ed s, de la forme:
\[(D_0\oplus \cdots\oplus D_0(n_1))\oplus (D_0(l)\oplus\cdots\oplus
D_0(l+n_2))\] avec $0\leq n_1\leq r_1$, $0\leq n_2\leq r_2$
vérifiant $n_1\leq l$ ou $l\leq n_1\leq l+n_2$.
\end{lemma}
\begin{proof}
Immédiat.
\end{proof}
Motivé par les exemples \ref{examples} et le lemme \ref{lemma-homo}
ci-dessus, on définit l'objet $(\varphi,N,\Gal(L'/L),D)$ de la
manière suivante:
\begin{itemize}
\item[--] en tant que $L_0'\otimes_{\Q_p}K$-module, $D=\oplus_{i=1}^sD_i$;

\item[--] on définit $N=\oplus_{i=1}^sN_i$ sur $D$;

\item[--] on définit l'action de $\Gal(L'/L)$ naturellement;

\item[--] pour tout $1\leq k_1<k_2\leq s$:
\begin{itemize}
\item[(a)] si $k_1,k_2\in H_i$ pour un $i$,  on peut écrire
\[D_{k_1}=D_0\oplus \cdots \oplus D_0(r_1),\ \ D_{k_2}=D_0(l)\oplus\cdots\oplus D_0(l+r_2);\]
avec $r_1,\ r_2,\ l\geq 0$. Si de plus $l=0$ ou $r_1=l+r_2$, et
$k_2$ est le plus petit entier de $H_i$ ayant l'une de ces
propriétés, on pose (sur $D_{k_1}\oplus D_{k_2}$)
\[\varphi (e_1)=\varphi_{k_1}(e_1)+\varphi_{k_2}(\alpha(e_1))\ \mathrm{si\ }e_1\in
D_{k_1}, \ \ \varphi (e_2)=e_2\ \mathrm{si\ }e_2\in D_{k_2}\] où\
$\alpha\neq 0$ est un élément fixé de
$\Hom_{\MOD_{L'/L}}(D_{k_1},D_{k_2})$.

\item[(b)] sinon, on pose $\varphi=\varphi_{k_1}\oplus
\varphi_{k_2}$ sur $D_{k_1}\oplus D_{k_2}$.
\end{itemize}
\end{itemize}
On déduit de la définition et du lemme \ref{lemma-homo} que
$(\varphi,N,\Gal(L'/L),D)\in \MOD_{L'/L}$ et
\[\WD(\varphi,N,\Gal(L'/L),D)^{\mathrm{ss}}=(r,N,V).\]

\subsection{Structure de $D$}\label{section-stru}
On conserve les notations du paragraphe précédent. Dans ce
paragraphe, on étudie la structure de $D$ qui sera utilisée dans
la suite.

D'abord, on suppose $D=D_{H_i}$ pour un $i$ de sorte qu'on peut
écrire
\[\begin{array}{rll}D&=&D_1\oplus\cdots \oplus D_s\\
&=&\big(D_0\oplus D_0(1)\oplus\cdots\oplus
D_0(b_1-1)\big)\oplus\cdots\oplus \big(D_0(l_s)\oplus\cdots\oplus
D_0(l_s+b_s-1)\big)
\end{array}\]  avec  $l_i\geq 0$ et $D_0$ absolument irréductible de dimension $h$.
%On suppose que $\dim D_0=1$.

%Par dé finition,
%$D=\bigoplus_{n=0}^{f'-1}V_{\sigma_0'\circ\varphi_0'^{-n}}$, donc on
%sait que $\varphi^{f'}$ est $L_0'\otimes_{\Q_p}K$-liné aire sur
%$D$.

 On considère pour l'instant
$(\varphi|_{D_1},N|_{D_1},\Gal(L'/L),D_1)\in\MOD_{L'/L}$. Par
hypothèse, $(\varphi|_{D_0},\Gal(L'/L),D_0)$ est un objet
absolument irréductible dans $\MOD_{L'/L}$, donc d'après la
remarque \ref{rem-scalaire}, $\varphi^{f'}|_{D_0}=a$ est scalaire
avec $a\in K^{\times}$. Notons $\varphi'=\varphi^{f'}$ et
$q'=p^{f'}$. %Alors $\varphi'$ est $L_0'\otimes_{\Q_p}K$-liné
%aire sur $D$.

On récrit $D$ sous la forme
\[\begin{array}{rll}D&=&D_{=0}\oplus\cdots\oplus D_{=n}\\
&=&\big(\bigoplus D_0\big)\oplus\big(\bigoplus
D_0(1)\big)\oplus\cdots\oplus(\bigoplus D_0(n)),\end{array}\] où
$n=\max_{2\leq i\leq s}\{r_1,l_i+b_i-1\}$. Alors par définition de
$\varphi$
\[D_{=j}=\{v\in D|\ (\varphi'-q'^{j}a)^kv=0, k\gg 0\}.\]

%\[D_{=k}^{N=0}=D_{=k}\cap \Ker N\]
%\[D_{=k}^{\im }=D_{=k}\cap \im N\]
\begin{defn}
(1) Un sous-objet $D'$ de $D$ est dit \emph{bon} s'il est de la
forme
\[\begin{array}{rll}D'&=&D'_1\oplus\cdots \oplus D'_s\\
&=&\big(D_0\oplus D_0(1)\oplus\cdots\oplus
D_0(b_1')\big)\oplus\cdots\oplus \big(D_0(l_s)\oplus\cdots\oplus
D_0(l_s+b_s')\big)
\end{array}\]
avec $-1\leq b_i'\leq b_i-1$, où $b_i'=-1$ désigne que $D'\cap
D_i=0$.

(2) On appelle \emph{drapeau} de $D$ une suite de
sous-$L_0\otimes_{\Q_p}K$-modules libres $\{E_i\}_{1\leq i\leq m}$
de $D$ tels que
\[0\subsetneq E_1\subsetneq\cdots\subsetneq E_m\subsetneq D;\]
%On  définit la notion de  \emph{sous-drapeau} d'un drapeau est de
%la manière évidente.
soit
\[0\subsetneq E'_1\subsetneq\cdots\subsetneq E'_{m'}\subsetneq D\]
un autre drapeau de $D$, il est dite un \emph{sous-drapeau} s'il
existe $1\leq j(i)\leq m'$ tel que $E_i=E_{j(i)}'$ pour tout
$1\leq i\leq m$.
 Soit $\sigma:L\hookrightarrow K$ un plongement, on définit
un drapeau de $D_{L',\sigma}$ de manière analogue.
%Un drapeau $\Delta$ est dit \emph{plus fin} que $\Delta'$ si

(3) Un drapeau
\[\Delta:0\subsetneq E_1\subsetneq\cdots\subsetneq E_m\subsetneq D\]
de $D$ est dit \emph{bon} si tous les $E_i$ sont bons.

%(iii) Un sous-objet de $D$ est dite \emph{bon} si elle est stable
%par $\varphi$, $N$ et $\Gal(L'/L)$, et son image sous $\pi:D\ra V$
%est bon.  Un drapeau de $D$ est dite \emph{bonne} si elle
%est stable par $\varphi$, $N$, $\Gal(L'/L)$, et son image sous
%$\pi:D\ra V$ est un bon drapeau  de $V$.
%(iii) Le drapeau
%\[0\subset D_{0,1}\subset\cdots\subset D_{0,?}\subset D_{1,1}\subset\cdots\subset D_{n,?}=D\]
\end{defn}
%D'aprè s le lemme, on fixe une filtration $\Fil$ sur $V$ telle que
%il est atapté e à\ tout les drapeaux bons.
\begin{rem}
(1) L'ensemble des bons sous-objets de $D$ est de cardinal fini.

(2) Si $E_1$ et $E_2$ sont deux bons sous-objets de $D$, alors
$E_1+E_2$, $E_1\cap E_2$, $\ker(N|_{E_1})$ et $N(E_1)$ sont bons
aussi. De plus, pour tout $j$, $\ker((\varphi'-q'^ja)|_{E_1})$ et
$(\varphi'-q'^ja)(E_1)$ sont bons.
\end{rem}

Si on définit
\[D_{=j,i}=\{v\in D_{=j}|\ (\varphi'-q'^{j}a)^iv=0\},\]
alors on obtient un bon drapeau de $D_{=j}$:
\[D_{=j,1}\subset D_{=j,2}\subset\cdots.\]
Posons $D_{<j}=\oplus_{i<j} D_{=i}$ et $D_{\leq j}=\oplus_{i\leq j}
D_{=i}$.

%\begin{rem}
%Si $D'$ est un bon sous-objet de $D$ et $a\in K^{\times}$, alors les
%objets $N(D')$ et $\ker N|_{D'}$, $(\varphi'-a) (D')$,
%$\ker(\varphi'-a)|_{D'}$ sont encore bons.
%\end{rem}

Dans le cas général $D=\oplus_{i}D_{H_i}$, on dit qu'un sous-objet
$D'$ est \emph{bon} s'il est somme directe de bons sous-objets de
$D_{H_i}$.

\subsection{Définition de la filtration sur $D$}
%On considé re un module filtré\ $D$ sur un anneau $R$ tel que
%tout sous-$R$-module de $D$ est libre.
%Soit $D$ un objet dans $\MOD_{L'/L}$ de rang $d+1$ (sur
%$L_0'\otimes_{\Q_p}K$).
%On fixe un plongement
%$\sigma:L\hookrightarrow K$ et considè re $D_{L',\sigma}$ comme un
%module sur $L'\otimes_{L,\sigma}K$.
On conserve les notations du paragraphe précédent. Soit
$\Fil_{\sigma}=(\Fil^iD_{L',\sigma})_{i}$ une filtration de
$D_{L',\sigma}$ stable par $\Gal(L'/L)$.  On note $I(\Fil_{\sigma})$
l'ensemble des entiers $i$ tels que
\[\Fil^{i}D_{L',\sigma}/\Fil^{i+1}D_{L',\sigma}\neq 0\]
avec multiplicité égale à $\dim_K
\Fil^{i}D_{L',\sigma}/\Fil^{i+1}D_{L',\sigma}$ (i.e. les opposés
des poids de Hodge-Tate).
\begin{rem}\label{rem-HT}
Si $D_1\subset D_2$ sont deux sous-$L'\otimes_{L,\sigma}K$-modules
de $D_{L',\sigma}$. Alors $I(\Fil_{\sigma},D_1)\\ \subset
I(\Fil_{\sigma},D_2)$.
\end{rem}
% Car
%\[\begin{array}{rll}&&i\in I(D_1) \mathrm{\ avec\ multiplicite\ }n\\
%&\Leftrightarrow& \dim D_1\cap\Fil^iD/D_1\cap\Fil^{i+1}D=n\\
%&\Rightarrow& \dim D_2\cap\Fil^iD/D_2\cap\Fil^{i+1}D\geq n\\
%&\Leftrightarrow&i\in I(D_2) \mathrm{\ avec\ multiplicite\ }\geq
%n.\end{array}\]

\begin{defn}\label{defn-adapte}
Soit \[\Delta_{\sigma}:0=E_{0,\sigma}\subsetneq
E_{1,\sigma}\subsetneq\cdots\subsetneq
E_{m+1,\sigma}=D_{L',\sigma}\] un drapeau de $D_{L',\sigma}$ par des
sous-$L'\otimes_{L,\sigma}K$-modules libres, et $\Fil_{\sigma}$ une
filtration sur $D_{L',\sigma}$ avec
$I(\Fil_{\sigma})=\{i_{1,\sigma}<i_{2,\sigma}<\cdots<i_{d+1,\sigma}\}$.
On dit que $\Fil_{\sigma}$ est \emph{transverse} à
$\Delta_{\sigma}$, si pour tout $k$,
\[I(\Fil_{\sigma},E_{k,\sigma})=\{i_{1,\sigma},\dots,i_{\dim_K{E_{k,\sigma}},\sigma}\}.\]
\end{defn}

\begin{lemma}\label{lemma-BS}
Soient  $\Delta=(\Delta_{\sigma})_{\sigma}$ un drapeau de $D$,  et
pour tout $\sigma:L\hookrightarrow K$ des entiers
$\{i_{j,\sigma}\}_{1\leq j\leq d+1}$ tels que $
i_{1,\sigma}<\cdots<i_{d+1,\sigma} $ . Alors il existe une
filtration $\Fil=(\Fil_{\sigma})_{\sigma}$ sur $D$ stable par
$\Gal(L'/L)$ telle que pour tout $\sigma$, $\Fil_{\sigma}$ est
transverse à $\Delta_{\sigma}$ et
\[I(\Fil_{\sigma})=\{i_{1,\sigma},\dots,i_{d+1,\sigma}\}.\]
De plus, si $K_1\subset K$ est un sous-corps de cardinal infini,
alors on peut prendre $\Fil_\sigma$ telle qu'elle soit définie sur
$K_1$.
\end{lemma}
\begin{proof}
Voir \cite{BS}, Proposition 3.2.
\end{proof}
\begin{rem}
De plus, si on se donne $\Delta_1$, $\dots$, $\Delta_n$  des
drapeaux de $D$ et des entiers
$\{i_{1,\sigma}<\cdots<i_{d+1,\sigma}\}$, on obtient encore par la
preuve du lemme \ref{lemma-BS} qu'il existe une filtration $\Fil$
telle que l'énoncé du lemme \ref{lemma-BS} soit vrai pour tout
$\Delta_i$.
\end{rem}

On fixe une filtration $\Fil=(\Fil_{\sigma})_{\sigma}$ sur $D$ qui
est transverse à tous les bons drapeaux (c'est un ensemble fini).
En particulier, comme tout bon sous-objet $D'$ de $D$ appartient à
un tel drapeau, on a (cf. la preuve de \cite{BS}, Proposition 5.1)
\[t_{H}(D'_{L'})=[K:L]\summ_{j=1}^{\dim D'}\summ_{\sigma}i_{j,\sigma}.\]

\begin{cor}\label{cor-decom-ineq}
Avec les notations précédentes.  Soient $\{E_i\}_{0\leq i\leq
m+1}$ des bons sous-objets de $D$ tels que $E_0=0$, $E_{m+1}=D$ et
$\dim E_{i-1}<\dim E_{i}$, et $D'$ un sous-objet quelconque de $D$.
Posons
\[a_{i }=\dim E_{i }-\dim E_{i-1 }, \ \ c_{i}=\dim (E_{i}\cap
D')-\dim (E_{i-1}\cap D'),\] et
\[\Omega=\{j\in\N |\mathrm{\ il\ existe\ }l  \mathrm{\ tel\ que\ }
\summ_{i=0}^{l}a_{i }-c_{l }+1\leq j\leq \summ_{i=0}^la_{i  }\}.\]
Alors
\[t_H(D'_{L'})\leq[K:L]\summ_{j\in \Omega}\summ_{\sigma}i_{j,\sigma}.\]
\end{cor}
\begin{proof}
Cela découle de la définition \ref{defn-adapte}, de la remarque
\ref{rem-HT} et du choix de $\Fil$.
\end{proof}
%\begin{lemma}
%Soient $D$ un $K$-espace vectoriel de dimension $d$, et $\Delta_1$,
%..., $\Delta_n$ des drapeaux de $D$. Alors il existe une
%filtration sur $D$ qui est transverse à\ tout les $\Delta_i$. De
%plus, si $K_1\subset K$ un sous-corps (donc a éléments
%infini), et $D=D'\otimes_{K'}K $, alors il existe une filtration
%$\Fil'$ sur $D'$ telle que $\Fil=\Fil'\otimes_{K'}K$ satisfait les
%conditions.
%\end{lemma}

\subsection{La preuve: cas spécial}
Dans ce paragraphe, on va prouver l'implication ``(ii)$\Rightarrow$
(iii)'' du théorème \ref{thm-main}
 %(sous l'hypothèse que $K$ est suffisamment gros)
 dans le cas spécial $D=D_{H_i}$. Rappelons qu'on a défini un certain objet
 $(\varphi,N,\Gal(L'/L),D)$ (cf. \S4.1) de $\MOD_{L'/L}$ tel que
 \[\WD(\varphi,N,\Gal(L'/L),D)^{\textrm{ss}}=(r,N,V),\]
 et on a fixé une filtration $\Fil=(\Fil^iD_{L',\sigma})_{i,\sigma}$ (cf. \S4.3) sur
 $D_{L'}$ telle que
 \[\Fil^iD_{L',\sigma}/\Fil^{i+1}D_{L',\sigma}\neq 0\Leftrightarrow i\in\{i_{1,\sigma},\dots,i_{d+1,\sigma}\}.\]
 Donc il suffit de prouver% cette filtration est \emph{admissible}.
 %Pour simplifier, on suppose $f'=1$, $L=\Q_p$ $D=D_{H_i}$ pour un $i$ et
%$\dim_K(D_0)=1$, alors $D_i=V_i$ en tant que $K$-espaces vectoriel
%grâce à\ $f'=1$.
\begin{theorem}\label{cor-admissible}
La filtration $\Fil$ sur $D_{L'}$ est admissible. Autrement dit, si
$D'$ est un sous-objet de $D$, alors
\[t_H(D'_{L'})\leq t_N(D').\]
\end{theorem}

On fixe un sous-objet $D'$ de $D$ et on commence par établir
certaines propriétés de $D'$.
\begin{defn}
Soient $E$, $E'$ deux bons sous-objets de $D$ tels que $\dim E'>\dim
E$, on définit
\[\alpha(E'/E,D')=\frac{\dim E'\cap D'-\dim
E\cap D'}{\dim E'-\dim E}.\]
 On note $\alpha(E'/E)=\alpha(E'/E,D')$ s'il n'y a pas de risque de confusion.
On écrit $\alpha(E')$ si $E=0$.
\end{defn}
\begin{rem}
(1) Le nombre $\alpha(E'/E)$ peut être négatif.

(2) Soient $E_1$, $E_2$, $E_3$ trois bons sous-objets de $D$ tels
que $\dim E_1<\dim E_2<\dim E_3$, alors ou bien $\alpha(E_3/E_2)>
\alpha(E_3/E_1)> \alpha(E_2/E_1)$, ou bien
$\alpha(E_2/E_1)>\alpha(E_3/E_1)> \alpha(E_3/E_2)$, ou bien
$\alpha(E_2/E_1)=\alpha(E_3/E_1)= \alpha(E_3/E_2)$.
\end{rem}
%On verrons prouver que pour un $(\varphi,N)$-sous-module $D'$ de
%$D$, on définit $D'_{=k,i}=D\cap D_{k,i}'$.

%\begin{lemma}\label{lemma-slope}
%Soient $\gamma:D\ra D$ un endomorphisme $L_0'\otimes_{\Q_p}K$-liné
%aire, et $D'$, $E$ deux sous-objets de $D$ tels que
%$\gamma(E)\subset E$. Posons $E^{\gamma=0}=\Ker(\gamma|_E)$. alors,
%on a
%\[0\ra E^{\gamma=0}\cap D'\ra E\cap D'\ra \gamma(E\cap D')\ra 0.\]
%En particulier, on a
%\[\dim E\cap D'\leq \dim E^{\gamma=0}\cap D'+\dim \gamma(E)\cap D'.\]
%\end{lemma}
%\begin{proof}
%Immédiat.
%\end{proof}
\begin{theorem}\label{thm-decom}
(1) Pour tout  sous-objet $D'$ of $D$, il existe un bon drapeau
\[\Delta=\Delta_{D'}: 0=E_0\subsetneq E_1\subsetneq\cdots\subsetneq E_{m+1}=D\]
 tel que

(a) pour $1\leq i\leq m$, $\alpha(E_{i+1}/E_i)\leq
\alpha(E_i/E_{i-1})$; si
$\alpha(E_{i+1}/E_i)=\alpha(E_i/E_{i-1})$, alors $\dim
E_{i+1}/E_i\geq\dim E_i/E_{i-1}$;

(b) $N(E_i)\subset E_{i-1}$;

(c) pour tout $i$, il existe un entier $j\leq 1$, tel que
$(\varphi'-q'^ja)E_i\subset E_{i-1}$;

(d) $\Delta$ n'a pas de sous-drapeau qui est bon et qui vérifie les
conditions (a)--(c).

(2)  Si
\[\Delta'=\Delta'_{D'}: 0=E'_0\subsetneq E'_1\subsetneq\cdots\subsetneq E'_{m'+1}=D\]
est un autre bon drapeau vérifiant les propriétés (a)--(d), alors
$m=m'$ et
\[\dim E_i=\dim E'_i, \ \
\alpha(E_{i+1}/E_{i})=\alpha(E'_{i+1}/E'_i)\]
 pour tout $i$.

%(v) pour tout $1\leq j\leq d+1$, on a $\dim V'\cap E$
%(v) Un tel drapeau est unique au sens que si
%\[\Delta': 0=E_0'\subset E_1'\subset\cdots\subset E_n'=V\]
%est une autre décompostion vérifiant les conditions (i)-(iii),
%alors il existe une bijection entre l'ensembles $\{E_{i+1}/E_i\}$ et
%$\{E_{i+1}'/E_i'\}$. Donc les pentes $\{\alpha_i\}$ ne dépend pas
%du choix de $\Delta$.
\end{theorem}

\begin{proof}
(1) On munit $\Q\times \N$ de l'ordre suivant: \[(x_1,y_1)\geq
(x_2,y_2)\Leftrightarrow (x_1> x_2) \textrm{ ou } (x_1=x_2 \textrm{
et } y_1\leq y_2).\]
 On considère
l'ensemble des bons sous-objets non nuls de $D$ et on choisit un
d'entre eux, noté $E_{1}$, tel que le couple $(\alpha(E_1), \dim
E_1)$ soit maximal. Puis on considère l'ensemble des bons
sous-objets de $D$ contenant $E_1$ strictement et on choisit un
d'entre eux, noté $E_2$, tel que $(\alpha(E_{2}/E_1),\dim E_2/E_1)$
soit maximal, et on obtient ainsi un bon drapeau
\[0=E_0\subsetneq E_1\subsetneq\cdots\subsetneq E_m\subsetneq E_{m+1}=D.\]

On prouve maintenant que  ce drapeau vérifie les conditions
(a)-(d):

%D'abord, on considè re l'ensemble des sous-objets bons et choisit
% un entre eux, $E_1$, tel que $\alpha_1(V')$ soit maximal. Aprè s, on
%choisit un objet bon $E_i$ contenant $E_1$ tel que $\alpha_2(V')$
%est maximal.

(a) et (d) sont automatiques et proviennent de la définition du
drapeau.

(b) Supposons d'abord $i=1$. Comme
\[\dim E_1=\dim N(E_1)+\dim E_1^{N=0}\]
et
\[ \dim (E_1\cap D')\leq  \dim N(E_1)\cap D'+\dim
(E_1^{N=0}\cap D'),\]% d'aprè s le lemme \ref{lemma-slope}.
 on déduit que ou bien $\alpha(N(E_i))\geq \alpha(E_1)$  ou bien $\alpha(E_1^{N=0})\geq
 \alpha(E_1)$, et donc,  par définition de
 $E_1$, $N(E_1)=0$ puisque $N(E_1)\neq E_1$.

 %Supposons que l'on a $N(E_{j})\subset E_{j-1}$ pour tout $j\leq i$, on prouverons que $N(E_{i+1})\subset E_{i}$.
Pour $i$ général, on considère $E_{i+1}'=\{v\in E_{i+1}|\ N(v)\in
E_{i}\}$ qui contient $E_i$. On a  une suite exacte de
$L_0'\otimes_{\Q_p}K$-modules libres
\[0\ra E_{i+1}'/E_i\ra E_{i+1}/E_i\ra (E_i+N(E_{i+1}))/E_i\ra 0,\]
donc d'après
%le lemme \ref{lemma-slope} et
le choix de $E_i$, on obtient
\[E_i+N(E_{i+1})=E_i \mathrm{\ ou\ } E_{i+1}.\]
Or $\ker (N|_{E_{i+1}})\neq 0$, on obtient $N(E_{i+1})\subseteq
E_i$.

(c) Supposons d'abord $i=1$. Comme
\[E_{1}=(E_1\cap D_{=0})\oplus\cdots \oplus (E_{1}\cap D_{=s}),\]
et
\[D'\cap E_1=(D'\cap (E_1\cap D_{=0}))\oplus\cdots\oplus (D'\cap (E_1\cap D_{=s})),\]
o\ud\ $D'\cap(E_1\cap D_{=j})$ est le sous-espace de $D'\cap E_1$
form\eb\ par les vecteurs $v$ tels que $(\varphi-q'^ja)^kv=0$ si
$k\gg 0$. Le choix de $E_1$ implique $E_1=E_1\cap D_{=j}$ pour un
$j$, i.e., $E_1\subset D_{=j}$. En utilisant la suite exacte
\[0\ra \ker((\varphi'-q'^{j}a)|_{E_1})\ra E_1\ra (\varphi'-q'^{j}a)({E_1})\ra 0,\]
on obtient $(\varphi'-q'^ja)(E_1)=0$. Dans le cas général, on
utilise l'argument analogue comme en (b).%, en notant que si on pose
%\[E_{i+1}'=\{v\in E_{i+1}| (\varphi'-q'^{j}a)v\in E_i\}\]
%et
%\[E_{i+1}''=E_i+(\varphi'-q'^ja)(E_{i+1}),\]
%alors $E_{i+1}'$ et $E_{i+1}''$ sont bons contenant $E_i$.

 (2) Cela sera une conséquence du lemme ci-dessous.
\end{proof}
\begin{lemma}\label{lemme-pour-decom}
On fixe un drapeau $\Delta_{D'}$ comme dans le théorème. Pour tout
bon sous-objet $L$ de $D$, soit $i\in\Z_{\geq0}$ l'unique entier tel
que
\[\dim E_i<\dim L\leq \dim E_{i+1}.\]
 Alors on a $\alpha(L/E_i)\leq \alpha(E_{i+1}/E_i)$.
 En particulier, si $\dim L=\dim E_{i+1}$, alors
 \[\dim D'\cap L\leq \dim D'\cap E_{i+1}.\]
\end{lemma}
\begin{proof}
 Supposons
l'énoncé faux et $L$ un tel sous-objet de $D$ de dimension
maximale. Soit $i$ l'unique entier tel que $\dim E_i<\dim L\leq \dim
E_{i+1}$. Alors par hypothèse, $\alpha(L/E_i)>\alpha(E_{i+1}/E_i)$.

Évidemment $L\neq D$. Soit $j$ le plus petit entier tel que
$E_{j+1}\nsubseteq L$, de sorte que $j\leq i$ et $E_j\subseteq L\cap
E_{j+1}$. On pose $L'=L+ E_{j+1}$ qui contient $L$ strictement et
soit $l$ l'unique entier tel que $\dim E_l<\dim L'\leq \dim E_{l+1}$
(alors $l\geq i$). On voit qu'il suffit de prouver que
\[\alpha(L'/E_l)>\alpha(E_{l+1}/E_{l}).\]

  Notons que $L\cap E_{j+1}$ est bon et contient $E_j$, on a
$\alpha (E_{j+1}/(L\cap E_{j+1}))\geq \alpha(E_{j+1}/E_j)$ par le
choix de $E_{j+1}$, et donc
\[\alpha(E_{j+1}/(L\cap E_{j+1}))\geq \alpha(E_{i+1}/E_i).\]%\ \ \ \forall k\geq j.\]
D'autre part, on a évidemment $\alpha(L'/L)\geq
\alpha(E_{j+1}/(L\cap E_{j+1}))$.  %alors $l\geq i$ et
% et $\alpha(L'/E_k)\leq \alpha(E_{k+1}/E_k)$. On a que $L_1$ satisfait que
Donc on obtient enfin $\alpha(L'/E_i)>\alpha(E_{i+1}/E_i)$, et ceci
achève la preuve dans le cas $l=i$. On est donc ramené au cas
$l>i$. Dans ce cas on a $\alpha(L'/E_i)>\alpha(E_{i+1}/E_i)\geq
\alpha(E_l/E_i)$, et ceci implique $\alpha(L'/E_l)>
\alpha(E_l/E_i)\geq \alpha(E_{l+1}/E_l)$.
\end{proof}
\begin{rem}
D'après (b) et (c) du théorème \ref{thm-decom}, on a pour tout
$i$, qu'il existe $j$, $l$ et une injection $E_{i+1}/E_i\ra
D_{=j,l+1}/D_{=j,l}$.
\end{rem}

On fixe un bon drapeau de $D$
\[\Delta=\Delta_{D'}:0=E_0\subsetneq E_1\subsetneq\cdots\subsetneq E_{m+1}=D\]
comme dans le théorème \ref{thm-decom}.

Soient $F_1$ le plus grand bon sous-objet tel que $\alpha(F_1)=1$ et
$F_2$ le plus grand bon sous-objet contenant $D'$. D'après le lemme
\ref{lemme-pour-decom}, on voit que $F_1$ et $F_2$ appartiennent à
$\Delta$. On note $k'$ et $k$ les  entiers tels que $F_1=E_{k'}$ et
$F_2=E_{k'+k}$ (indépendants du choix du drapeau). On considère le
drapeau
\[0=E_0\subsetneq E_{k'}=F_1\subsetneq\cdots\subsetneq E_{k'+k}=F_2\subsetneq E_{m+1}=D,\]
et on pose $a_0=\dim F_1$, $a_{k+1}=\dim D-\dim F_2$, et pour $1\leq
i\leq k$
 \[ a_i=\dim E_{k'+i}/E_{k'+i-1},\ \ \  c_i=\dim D'\cap E_{k'+i}/D'\cap
E_{k'+i-1}.\]
\begin{prop}\label{thm-special}
Avec les notations plus haut, on a
\[ a_0\geq \max\{c_i\}_{1\leq i\leq k},\ \ \   a_{k+1}\geq \max\{a_i-c_i\}_{1\leq i\leq k}.\]
\end{prop}
\begin{proof}
On prouve seulement la première inégalité, la preuve de l'autre
étant analogue. D'aprè s l'assertion (i) du lemme
\ref{lemma-compare-dim} ci-dessous, on peut supposer $h=\dim
D_0=1$.

 Choisissons un $r$  tel que
$c_{r}=\max_{1\leq j\leq k}\{c_j\}$. On sait qu'il existe $j$, $l$
tels que $ E_{k'+r}/E_{k'+r-1}\hookrightarrow D_{=j,l+1}/D_{=j,l}$.
On a donc $\dim D'\cap D_{=j,l+1}-\dim D'\cap D_{=j,l}\geq c_r$, et
puis $\dim D'\cap D_{=j,1}\geq c_r$ parce que l'on a un
homomorphisme injectif
\[(\varphi'-q'^ja): D_{=j,l+1}/D_{=j,l}\ra D_{=j,l}/D_{=j,l-1}.\]
D'après le lemme \ref{lemma-compare-dim} ci-dessous, on a pour tout
$1\leq i\leq j$,
\[\dim D'\cap D_{=i-1,1}\geq \dim D'\cap D_{=i,1}-1,\]
 et pour que $\dim D'\cap D_{=i-1,1}=\dim D'\cap D_{=i,1}-1$, il faut
 $\dim D_{=i-1,1}=\dim D_{=i,1}-1$, d'où
 \[\dim \ker(N|_{D_{=i,1}})=1,\ \ \ \ker
N|_{D_{=i,1}}\subset D'.\] S'il en est ainsi, alors par la
définition de $F_1$ on a
\[\ker( N|_{D_{=i,1}})\subset F_1,\]
et ceci permet de conclure  puisque $\dim D'\cap D_{=1,1}\leq 1$
d'après le lemme \ref{lemma-compare-dim}.
%on prouve que si $a_i-b_i$ est
%maximal, et $E_{i}/E_{i-1}\hookrightarrow D_{=j,l+1}/D_{=j,l}$ pour
%$j$, $l$ convenables, alors $D_{=h,l+1}\neq 0$ et
%\[\dim D_{=h,l+1}/D_{=h,l}-\dim D'\cap D_{=h,l+1}/D'\cap D_{=h,l}\geq a_i-b_i.\]
\end{proof}

\begin{lemma}\label{lemma-compare-dim}
(1) Pour tout sous-objet $\overline{D}$ de $D$, on a $h|\dim
\overline{D}$.

(2) Pour tout $k$, $\dim \ker N\cap \ker (\varphi'-q'^ka)=h$ ou $0$.

(3) $\dim D_{=1,1}=h$, $|\dim D_{=i,1}-\dim D_{=i-1,1}|= h$ ou $0$.
\end{lemma}
\begin{proof}
(1) Comme $\overline{D}=\bigoplus_{j=1}^s \overline{D}\cap
D_{=j}$, on peut supposer $N=0$ et alors par définition,
$(\varphi,\Gal(L'/L),D)$ est une extension successive de objets de
type $(\varphi|_{D_0},\Gal(L'/L),\\ D_0)$ qui est absolument
irréductible. L'énoncé s'en déduit.

(2) Cela résulte de la définition de $\varphi$ et de (1).

(3) Le premier énoncé est une conséquence directe de (2). Pour le
deuxième, on considère $N(D_{=i,1})$. D'une part, on a $\dim
D_{=i,1}-\dim N(D_{=i,1})\leq h$ d'après (2). D'autre part, $\dim
D_{=i-1,1}-\dim N(D_{=i,1})\leq h$ d'après la définition de
$\varphi$. Enfin,  l'assertion (i) nous permet de conclure.
\end{proof}

 En vue du théorème \ref{thm-decom} et de la proposition \ref{thm-special}, on a besoin de la définition suivante:

\begin{defn}\label{def-special}
Soient $\alpha=(a_0,\dots,a_{k+1})$ une partition de $d+1$ et
$\Omega$ un sous-ensemble de $\{1,\dots,d+1\}$. Posons
\[I_0=\{1,2,\dots,a_0\}, I_1=\{a_0+1,\dots,a_0+a_1\},\cdots, I_{k+1}=\{a_0+\cdots+a_k+1,\dots,d+1\},\]
et  $c_i=|\Omega\cap I_i|$  pour $0\leq i\leq k+1$. On dit que la
paire $(\Omega,\alpha)$ est de type \emph{spécial} si les nombres
$\{a_i\}_{0\leq i\leq k+1}$ et $\{c_i\}_{0\leq i\leq k+1}$
vérifient les conditions suivantes:

(i) $a_0=c_0>0,\ c_{k+1}=0$, et $a_i\geq c_i> 0$ pour $1\leq i \leq
k$,

(ii)  $\frac{c_i}{a_i}\geq \frac{ c_{i+1}}{a_{i+1}}$ pour tout $i$,

(iii) $a_0\geq \max_{1\leq i\leq k}\{c_i\}$, $a_{k+1}\geq
\max_{1\leq i\leq k}\{a_i-c_i\}$.
\end{defn}

Les paires de types spéciaux ont la propriété suivante:
\begin{prop}\label{prop-special}
Soient $(\Omega,\alpha)$ une paire de type spécial, et
$\{m_i\}_{1\leq i\leq d+1}$, $\{n_i\}_{1\leq i\leq d+1}$ des nombres
réels tels que

(a) $m_i\leq m_{i+1}$, $n_i\leq n_{i+1}$,

(b) $m_{i+1}-m_i\geq n_{i+1}-n_i$,

(c) $\summ_{i=1}^{d+1}m_i\leq \summ_{i=1}^{d+1}n_i$.

Alors on a $\summ_{i\in\Omega}m_i\leq\summ_{i\in\Omega}n_i$.
\end{prop}
\begin{rem}
Cette proposition peut être vue comme une généralisation de
\cite{BS}, lemme 5.4.
\end{rem}

La démonstration du théorème \ref{cor-admissible} utilise la
proposition \ref{prop-special}, dont la preuve sera donnée au
 paragraphe suivant.

\begin{proof}[Démonstration de \ref{cor-admissible}]
Reprenons les notations de la proposition \ref{thm-special}. Donc on
a $\summ_{i=0}^{k+1} a_i=d+1$. Notons $\alpha=(a_0,\dots,a_{k+1})$
la partition de $d+1$,
\[I_0=\{1,\dots,a_0\},\cdots, I_{k+1}=\{a_0+\cdots+a_k+1,\dots,d+1\},\]
et $\Omega$ le sous-ensemble de $\{1,\dots,d+1\}$ tel que
\[ \Omega\cap I_i=\{\summ_{j=0}^{i}a_j-c_i+1,\summ_{j=0}^{i}a_j-c_i+2,\dots,\summ_{j=0}^ia_j\}.\]
On sait que, d'après le théorème \ref{thm-decom} et la
proposition \ref{thm-special}, $(\Omega,\alpha)$ est de type
spécial.

On pose $s_i=\dim (D'\cap D_{= i})$, $r_i=\sharp\{j\in\Omega|\ \dim
D_{<i}+1\leq j\leq \dim D_{\leq i}\}$ et
\[n_{j}=\left\{ {\begin{array}{ll} t_N(D_0)& 1\leq j\leq\dim{D_{=0}}\\
t_N(D_0(1)) &\dim D_{=0}+1\leq j\leq\dim D_{\leq 1}\\
\vdots&\vdots\\
t_N(D_0(n))&\dim D_{<n}+1\leq j\leq \dim D_{\leq n}=d+1 \\
\end{array}}\right.\]
de sorte qu'on a (notons que $t_N(D_0(i))=t_N(D_0)+i[K:\Q_p]$)
\[t_N(D')=\summ_{i=0}^ns_in_i=\dim D' \cdot t_N(D_0)+[K:\Q_p]\summ_{i=0}^nis_i,\]
et \[\summ_{j\in\Omega}n_j=\summ_{i=0}^nr_in_i=\dim D'\cdot
t_N(D_0)+[K:\Q_p]\summ_{i=0}^nir_i.\]

(a) D'abord, on prouve le théorème sous l'hypothèse
$\summ_{i=0}^ls_i\leq \summ_{i=0}^lr_i$ pour tout $l$. Alors on a:
\begin{itemize}
\item[--] $t_H(D'_{L'})\leq [K:L]\summ_{j\in
\Omega}\summ_{\sigma}i_{j,\sigma}$. Ceci est une conséquence du
corollaire \ref{cor-decom-ineq}.

\item[--]
$[K:L]\summ_{j\in \Omega}\summ_{\sigma}i_{j,\sigma}\leq \summ_{j\in
\Omega}n_j$. En posant $m_j=[K:L]\summ_{\sigma}i_{j,\sigma}$, on
voit que \[m_{j+1}-m_j\geq [K:\Q_p]\geq n_{j+1}-n_j\] puisque que
 $\dim D_{=i}\neq 0$ pour tout $0\leq i\leq n$. On déduit l'énoncé de la proposition
\ref{prop-special}.

\item[--] $t_N(D')\geq \summ_{j\in \Omega}n_j$. Il suffit de prouver que
$\sum_{i=0}^ni(s_i-r_i)\geq 0$. Mais
\[\begin{array}{rll}\summ_{i=0}^ni(s_i-r_i)&=&n\summ_{i=0}^n(s_i-r_i)-\summ_{i=0}^n(n-i)(s_i-r_i)\\
&=&n\summ_{i=0}^n(s_i-r_i)-\summ_{l=0}^{n-1}\summ_{i=0}^l(s_i-r_i),\end{array}\]
l'inégalité résulte du fait que
$\summ_{i=0}^ns_i=\summ_{i=0}^nr_i$ et $\summ_{i=0}^l(s_i-r_i)\leq
0$ par hypothèse.
\end{itemize}
Donc on obtient $t_H(D'_{L'})\leq t_N(D')$.

%\[\Omega=\{j\in \N|\ \textrm{il existe un } l \textrm{ tel que }\dim E_{l}-b_l+1\leq j\leq \dim E_{l}\},\]
(b) Dans le cas général,  on voit facilement qu'il suffit de
trouver un autre sous-ensemble $\Omega'$ de $\{1,\dots,d+1\}$ tel
que
\begin{itemize}
\item[--] $t_H(D'_{L'})\leq
[K:L]\summ_{j\in\Omega'}\summ_{\sigma}i_{j,\sigma}$;

\item[--] $(\Omega',\alpha)$ est de type spécial et $|\Omega'\cap I_i|=|\Omega\cap
I_i|=c_i$;

\item[--] pour tout $l$
\[\dim(D'\cap D_{\leq l})\leq
\sharp\{j\in\Omega'|j\leq \dim D_{\leq l}\}.\leqno(**)\]
\end{itemize}

 Si $0\leq m\leq n$, on notera $i(m)$ l'unique entier tel
que $0< i(m)\leq k$ et $\dim E_{i(m)-1}<\dim D_{\leq m}\leq \dim
E_{i(m)}$. D'abord, soit $m_1$ le plus petit entier tel que
\[\dim(D'\cap D_{\leq m_1})>\sharp\{j\in\Omega|j\leq \dim D_{\leq m_1}\}.\]
Alors si on pose
\[x_{m_1}=\dim D'\cap D_{\leq m_1}-\dim D'\cap E_{i(m_1)-1}\]
et \[ y_{m_1}=(\dim E_{i(m_1)}-\dim D_{\leq m_1})-(\dim D'\cap
E_{i(m_1)}-\dim D'\cap D_{\leq m_1}),\] on en déduit que $x_m>0$ et
$y_m>0$.
 On définit un sous-ensemble $\Omega'_1$ de $\{1,\dots,d+1\}$,
en modifiant la définition de $\Omega$, tel que
 $\Omega_1'\cap I_j=\Omega\cap I_j$ pour
$j\neq i(m_1)$ et
%\[\begin{array}{rll}\Omega'\cap[\dim E_{i({m_1})-1}+1,\dim E_{i(m_1)}]&=&\{j\in\N|\ \dim D_{\leq{m_1}}-x_{m_1}+1\leq j\leq \dim D_{\leq
%m_1}, \textrm{ou }\\
%&&\dim E_{i(m_1)}-(b_{i(m_1)}-x_{m_1})+1\leq j\leq \dim E_{i(m_1)}
%\}.\end{array}\]
\[\begin{array}{rll} \Omega_1' \cap I_{i(m_1)}&=&\{\dim D_{\leq{m_1}}-x_{m_1}+1,\dots,\dim D_{\leq m_1}\}\sqcup\\
&&\{\dim E_{i(m_1)}-(c_{i(m_1)}-x_{m_1})+1,\dots,\dim
E_{i(m_1)}\}.\end{array}\]
 Notons que $c_{i(m_1)}-x_{m_1}> 0$ d'après
l'hypothèse sur $m_1$ et le théorème \ref{thm-decom}, (2). Après,
si $\Omega'_1$ vérifie $(**)$, on pose $\Omega'=\Omega_1'$; sinon,
soit $m_2$ le plus petit entier tel que $(**)$ soit faux pour
$\Omega_1'$, alors $m_2>m_1$ et $y_{m_2}>0$ où
\[ y_{m_2}=(\dim E_{i(m_2)}-\dim D_{\leq m_2})-(\dim D'\cap E_{i(m_2)}-\dim D'\cap D_{\leq m_2}).\]
 Il y a deux possibilités:
\begin{itemize}
\item[--] si $i(m_2)>i(m_1)$, on pose
\[x_{m_2}=\dim D'\cap D_{\leq m_2}-\dim D'\cap E_{i(m_2)-1},\]
l'hypothèse sur $m_2$ implique que $x_{m_2}>0$;  on définit $\Omega_2'$ comme précédemment; %(mais $\Omega_2'\cap
%I_{j}=\Omega_1'\cap I_j$ si $j\neq i(m_2)$);

\item[--] si $i(m_1)=i(m_2)$, on pose
\[x_{m_2}=\dim D'\cap D_{\leq m_2}-\dim D'\cap D_{\leq m_1},\]
de sort que $x_{m_2}>0$;  on définit $\Omega_2'$ en modifiant la
définition de $\Omega_1'$ par
\[\begin{array}{rll}&&\Omega'_2\cap\{\dim D_{\leq m_1}+1,\dots,\dim E_{i(m_1)}\}\\
&=&\{\dim D_{\leq{m_2}}-x_{m_2}+1,\dots,\dim D_{\leq
m_2}\}\sqcup\\
&&\{\dim E_{i(m_1)}-(c_{i(m_1)}-x_{m_1}-x_{m_2})+1,\dots,\dim
E_{i(m_1)}\}.\end{array}\]
\end{itemize}
 Par récurrence on obtient bien un sous-ensemble $\Omega'$ de
 $\{1,\dots,d+1\}$ vérifiant les conditions
 demandées: les deux dernières  résultent de la définition et la première résulte encore
 du corollaire \ref{cor-decom-ineq} puisque
 $D_{\leq i}$ est aussi un bon sous-objet de $D$.

%On dé finit l'ensemble $E'\subset\{1,\cdots,d+1\}$ tel que
%\[\begin{array}{ll}j\in E'\cap [\dim E_{i-1}+1,\dim E_{i}] \Longleftrightarrow\\
% \left\{ {\begin{array}{ll}
%\dim E_{i}-b_i+1\leq j\leq \dim E_{i}& \mathrm{si\ }i\neq i(m)  \\
%\dim E_{i}-b_i+1\leq j\leq \dim E_{i}& \mathrm{si\ } i=i(m), x_m\leq y_m \textrm{\ où\ }x_m\leq 0\\

% {\begin{array}{l}\dim D_{\leq m}-x_m+1\leq j\leq \dim D_{\leq m}\textrm{\ o\`{u}} \\
% \dim E_{i(m)}-(b_{i(m)}-x_m)+1\leq j\leq \dim E_{i(m)}\end{array}}
% &\mathrm{si\ }i=i(m) \mathrm{\ et\ } x_m>\max\{y_m,0\}.\\
%\end{array}}\right.\end{array}\]
%(Notons que $b_{i(m)}-x_m> 0$ dans le derniè r  cas).

 \end{proof}

\subsection{Preuve de \ref{prop-special} }\label{un lemme}
Dans ce paragraphe, on prouve la proposition \ref{prop-special}. On
commence par un lemme:

\begin{lemma}\label{lemma-t}
Soient $k\geq1$  et $\{a_i\}_{0\leq i\leq k+1}$, $\{c_i\}_{0\leq
i\leq k+1}$ des nombres réels vérifiant les conditions (i)-(iii)
de la définition \ref{def-special}, alors  il existe des nombres
réels $\{t_i\}_{1\leq i\leq k}$ tels
 que

(i)'
$\frac{t_1}{a_0}=\frac{t_2}{c_1}=\cdots=\frac{t_k}{c_{k-1}}=:r$,

(ii)' pour tout $1\leq l\leq k$,
\[\summ_{i=1}^{l}{t_i}\geq \summ_{i=1}^l(a_i-c_i),\]

(iii)' %$\frac{t_1}{a_1}\leq r$, et
$\summ_{i=1}^{k}{t_i}-\summ_{i=1}^k(a_i-c_i)+rc_{k}\leq a_{k+1}$.
\end{lemma}

\begin{proof}
On peut écrire, si les $t_i$ existent, $t_i=t_1c_{i-1}/a_0$ pour
$1\leq i\leq k$, et alors (ii)' se réduit à
\[ t_1(1+\frac{1}{a_0}\summ_{i=1}^{l-1}c_i )\geq \summ_{i=1}^{l}(a_i-c_i),\]
et (iii)' se réduit à
\[%x\leq a_1,\ \ \ \
t_1(1+\frac{1}{a_0}\summ_{i=1}^{k}c_i )\leq
a_{k+1}+\summ_{i=1}^k(a_i-c_i).\]

On prend $t_1=\max_{1\leq l\leq
k}\{(\summ_{i=1}^{l}(a_i-c_i))(1+\frac{1}{a_0}\summ_{i=1}^{l-1}c_i
)^{-1}\}$, et $t_i=t_1c_{i-1}/a_0$, alors (i)', (ii)' sont
satisfaits. Pour voir que les $t_i$ ainsi définis satisfont à
(iii)', il suffit de vérifier que pour $1\leq l\leq k$,
\[(\summ_{i=1}^{l}(a_i-c_i))(1+\frac{1}{a_0}\summ_{i=1}^{k}c_i )\leq (1+\frac{1}{a_0}\summ_{i=1}^{l-1}c_i
) (a_{k+1}+\summ_{i=1}^k(a_i-c_i)),\] qui est une conséquence
facile du (ii) et (iii).
\end{proof}
Dans la suite de cet article, si on se donne les $\{a_i\}_{0\leq
i\leq k+1}$ et $\{c_i\}_{0\leq i\leq k+1}$ comme dans le lemme
\ref{lemma-t}, on prendra les $\{t_i\}_{1\leq i\leq k}$ et $r$ comme
dans le lemme tel que $r$ est le plus petit possible.

%\begin{proof}
%Le lemme ci-dessus nous permet de définir des nombres $(t_i)_{1\leq
%i\leq k}$ vérifiant les conditions (i)'-(iii)'. Posons
%\[I_i' =\left\{ {\begin{array}{ll}
%[a,a+a_0[& \mathrm{si\ }i=0 \\
%
% [a+a_0+\summ_{j=1}^{i-1}(t_i+b_j),a+\summ_{j=1}^{i}(t_j+b_j)[ &\mathrm{si\ }1\leq i\leq k\\
%
%[a+\summ_{j=1}^{k}(t_j+b_j),a+\summ_{j=1}^{k}(t_j+b_j)+rb_k[&\mathrm{si\
%}i=k+1,\end{array}}\right.\] et
%\[J_i'=\left\{ {\begin{array}{ll}
%[a,a+a_0[& \mathrm{si\ }i=0 \\
%
%\summ_{j=1}^i(t_j-(a_j-c_j))+[(a+\summ_{j=0}^{i}a_j)-c_i,a+\summ_{j=0}^{i}a_j[ &\mathrm{si\ }1\leq i\leq k\\
%
%\emptyset&\mathrm{si\ }i=k+1.\end{array}}\right.\] On a alors
%$|J_i|=|J_i'|$.
% Les conditions (ii)' et (iii)' impliquent que
%
%(a) $J_i\subset J_i'\subset I_i'\subset I$,
%
%(b) pour $1\leq i\leq k$, $|J_i'|=b_i$, $|I_i'|=t_i+b_i$, et
%$|I_{k+1}'|=rb_k$.\\
%Posons $I'=\bigsqcup_{i=0}^{k+1} I_i'\subset I$.
%
% Par hypothè se, $ \int_{I}(f-g)\leq 0$, on obtient que $ \int_{I'}(f-g)\leq
% 0$ puisque $f'\geq g'$, et donc
% \[\begin{array}{rll}0\geq & \summ_{j=0}^{k}\int_{I_i'}(f-g) \\
% = & \summ_{j=0}^{k}\int_{J_i'}(f-g)+\summ_{j=1}^{k+1}\int_{I_i'\backslash J_i'}(f-g)\\
% \geq & (1+r)\summ_{j=0}^k\int_{J_i'}(f-g),\end{array}\]
%ce qui implique que $\summ_{j=0}^k\int_{J_i'}(f-g)\leq 0$. On a donc
%\[ \summ_{j=0}^k\int_{J_i'}(f-g)\leq 0,\] et ensuite
%%$f'-g'\geq 0$, on obtient
%\[ \summ_{j=0}^k\int_{J_i}(f-g)\leq \summ_{j=0}^k\int_{J_i'}(f-g)\leq 0.\]
%Les conditions (i)'-(iii)' impliquent que
%\[\summ_{j=0}^k\int_{J_i}(f-g)\leq \summ_{j=0}^k\int_{J_i'}(f-g)\leq \frac{1}{1+r}\int_I(f-g)\leq 0.\]
%\end{proof}

\begin{proof}[Démonstration de \ref{prop-special}]
Grâce à la condition (b): $m_{i+1}-m_i\geq n_{i+1}-n_i$, on voit
qu'on peut supposer
\[\Omega=\{1\leq j\leq d+1|\ \textrm{il existe }l\textrm{ tel que }\summ_{i=0}^la_i-c_l+1\leq j\leq \summ_{i=0}^la_i\}.\]
Posons $I=[0,d+1[$. On prend deux fonctions lisses $f,g:I\ra\R$
telles que (c'est toujours possible)
\begin{itemize}

\item[--] $f'(x)\geq g'(x)\geq 0$, et $\int_{I}f\leq \int_{I}g$, et

\item[--] pour $0\leq i\leq d$, $f(i)=m_{i+1}$, $g(i)=n_{i+1}$.
\end{itemize}
Le lemme ci-dessus nous permet de définir des nombres $(t_i)_{1\leq
i\leq k}$ vérifiant les conditions (i)'-(iii)'. Posons
\[J_i =\left\{ {\begin{array}{ll}
[0,a_0[& \mathrm{si\ }i=0 \\

[\summ_{j=0}^{i}a_j-c_i,\summ_{j=0}^{i}a_j[ &\mathrm{si\ }1\leq i\leq k\\

\emptyset&\mathrm{si\ }i=k+1,\end{array}}\right.\] et
%\[I_i' =\left\{ {\begin{array}{ll}
%[a,a+a_0[& \mathrm{si\ }i=0 \\
%
% [a+a_0+\summ_{j=1}^{i-1}(t_i+b_j),a+\summ_{j=1}^{i}(t_j+b_j)[ &\mathrm{si\ }1\leq i\leq k\\
%
%[a+\summ_{j=1}^{k}(t_j+b_j),a+\summ_{j=1}^{k}(t_j+b_j)+rb_k[&\mathrm{si\
%}i=k+1,\end{array}}\right.\] et
\[J_i'=\left\{ {\begin{array}{ll}
[0,a_0[& \mathrm{si\ }i=0 \\

[a_0+\summ_{j=1}^{i-1}(t_j+c_j)+t_i,a_0+\summ_{j=0}^{i}(a_j+c_j)[ &\mathrm{si\ }1\leq i\leq k\\

\emptyset&\mathrm{si\ }i=k+1.\end{array}}\right.\] Alors on a
$|J_i|=|J_i'|$, $\bigsqcup_{1\leq i\leq k } J_i'\subset I$ et
\[\summ_{i\in \Omega}(m_{i}-n_i)\leq \summ_{i=0}^k\int_{J_i}(f-g)\leq \summ_{i=0}^k\int_{J_i'}(f-g)
\leq \frac{1}{r+1}\int_{I}(f-g)\leq 0,\] ce qui permet de conclure.
%\[ \summ_{i=0}^k\int_{J_i}(f-g)\leq 0.\]
%Mais on voit que
% \[\begin{array}{rll}\summ_{i=0}^k\int_{J_i}(f-g)&\geq &\summ_{i=0}^k\summ_{j=1}^{b_i}(f-g)(x_j')\\
% &\geq&\summ_{i=0}^k\summ_{j=1}^{b_i}(f-g)(\alpha_i+j-1)\\
% &=&\summ_{i\in \Omega}(m_{i}-n_i)\end{array} \]
% avec $x_j'\in [\alpha_i+j-1,\alpha_i+j[$, où\  on a é crit $\alpha_i=\summ_{j=0}^ia_j-b_i$ pour
% simplifier. Ceci permet de conclure.
\end{proof}
\subsection{Le cas général}
On considère le cas général: $D=\bigoplus_{i=1}^v D_{H_i}$.
%On commence par établir une généralisation de la
D'abord, la proposition \ref{prop-special} se généralise
aisément:
\begin{prop}
Soit $\alpha=(d_1,\dots,d_v)$ une partition de $d+1$, et pour tout
$1\leq i\leq v$ soient
$\alpha_i$ une partition de $d_i$ % où\
%\[I_i=\{\summ_{j=1}^{i-1}d_{j}+1,\dots,\summ_{j=1}^id_j\}.\]
 et $\Omega_i$ un sous-ensemble de $\{1,\dots,d_i\}$ tels que la
 paire
$(\Omega_i,\alpha_i)$ est de type spécial. On associe, à chaque
paire $(\Omega_i,\alpha_i)$, le nombre réel positif $r_i$ comme
dans le lemme \ref{lemma-t}, et on suppose $r_1\geq r_2\geq
\cdots\geq r_v$. Posons
\[\Omega_i'=\{d_1+\cdots d_{i-1}+l|\ l\in\Omega_i\},\
 \  \ \Omega=\bigsqcup_i \Omega_i'\] de sorte que $\Omega$
est un sous-ensemble de $\{1,\dots,d+1\}$. De plus, soient
$\{m_i\}_{1\leq i\leq d+1}$ et $\{n_i\}_{1\leq i\leq k+1}$ des
nombres réels vérifiant les conditions (a)-(c) de la proposition
\ref{prop-special}, alors on a
$\summ_{j\in\Omega}m_j\leq\summ_{j\in\Omega}n_j$.
\end{prop}
La preuve de cette proposition est analogue à celle de
\ref{prop-special} et on laisse les détails au lecteur.

Maintenant, posons $d_i=\dim D_{H_i}$, et $D_i'=D'\cap D_{H_i}$
pour tout $1\leq i\leq v$ et tout sous-objet $D'$ de $D$. Comme
$D_i'$ est un sous-objet de $D_{H_i}$, on peut lui associer un
drapeau $\Delta_i=\Delta_{D_i'}$ d'après le théorème
\ref{thm-decom}, et puis un sous-ensemble $\Omega_i$ de
$\{1,\dots,d_i\}$ et un nombre réel positif $r_i$ comme dans le
lemme \ref{lemma-t}. On peut supposer $r_1\geq r_2\geq\cdots\geq
r_v$. Posons
\[\Omega_i'=\{d_1+\cdots d_{i-1}+l|\ l\in\Omega_i\},\
\textrm{ et } \ \Omega=\bigsqcup_i \Omega_i'.\] On a alors
$t_H(D')\leq \summ_{j\in\Omega}\summ_{\sigma}i_{j,\sigma}$ par le
corollaire \ref{cor-decom-ineq}. Un argument analogue à la preuve
du corollaire \ref{cor-admissible}, en utilisant la proposition
ci-dessus, nous permet de conclure.

%Finalement, bien que on a supposé\ que
%$\dim_{L_0'\otimes_{\Q_p}K}D_0=1$, l'argument dans \S4.4 est valable
%pour le cas général.
%\begin{cor}
 %Pour tout $1\leq j\leq s$,  soient $(m_{ji})_{1\leq i\leq d_j}$,
 %$(n_{ji})_{1\leq i\leq d_j}$ des nombres réels vérifiant les
 %conditions (i)-(iii) de la proposition \ref{prop-special}, et
%$(a_{ji})$
%$(b_{ji})$ des nombres entiers positifs comme dans la proposition
%ci-dessus. On pose
%\[\Omega_j=\{i\in \N|\ \textrm{il existe }l,\ \textrm{tel que }\lambda_j+\summ_{r=0}^la_{jr}-b_{jl}+1\leq i\leq\lambda_j+\summ_{r=0}^{l}a_{jr}\}\]
%et $\Omega=\sqcup\Omega_j$. Alors
%\[\summ_{j=1}^s\summ_{i\in \Omega_j}m_{ji}\leq\summ_{j=1}^s\summ_{i\in\Omega_j}n_{ji}.\]
%\end{cor}
%\begin{proof}
%La mê me preuve que cette de la proposition \ref{prop-special}.
%\end{proof}

%\renewcommand{\refname}{Références}

 \hspace{5cm} \hrulefill\hspace{5.5cm}

Département de Mathématiques, bâtiment 425,  Université de
Paris-Sud, 91405 Orsay cedex, France

E-mail: yongquan.hu@math.u-psud.fr

\end{document}